\documentclass[smallcondensed,envcountsect]{svjour3}
\usepackage{geometry}                
\usepackage{graphicx}
\usepackage{amssymb}
\usepackage{epstopdf}
\usepackage{algorithm2e}
\usepackage[reqno]{amsmath}
\usepackage{amssymb}
\usepackage[]{graphicx}               
\usepackage{bm}
\usepackage{esvect}
\usepackage{ulem}
\usepackage{adjustbox}
\usepackage{accents}
\usepackage{stmaryrd}
\usepackage{xcolor}
\usepackage{enumitem}
\usepackage{hyperref}
\usepackage{booktabs}
\usepackage{longtable}
\usepackage[section]{placeins}
\usepackage{tabu}
\usepackage{pgfplots}

\usepackage{color}
\definecolor{chcol}{rgb}{0.4,0.,0.9}

\newcommand{\halfComma}{\kern 0.083334em}

\newcommand{\jump}[1]{\ensuremath{\left\llbracket #1 \right\rrbracket}}
\newcommand{\dS}{{\,\operatorname{dS}}}         
\newcommand{\PBT}{{\,\operatorname{PBT}}}    
\newcommand\iprod[1]{\left\langle #1\right\rangle}                                             
\newcommand\inorm[1]{\left |\left| #1\right|\right|}                                               
\newcommand\spacevec[1]{\accentset{\,\rightarrow}{#1}}                        
\newcommand\statevec[1]{\mathbf #1}                                                     
\newcommand\statevecGreek[1]{\boldsymbol #1}                                     
\newcommand\acclrvec[1]{\accentset{\,\leftrightarrow}{#1}}                      
\newcommand\blockvec[1]{\acclrvec{{\mathbf #1}}}                                  
\newcommand\mmatrix[1]{\underbar{#1}}				
\newcommand\bigmatrix[1]{\mathfrak #1}                          

\newcommand\oneHalf{\frac{1}{2}}

\newcommand\sym[1]{ #1^{s}}                                             

\usepackage{xargs}    
\usepackage[colorinlistoftodos,prependcaption,textsize=tiny]{todonotes}
\newcommandx{\unsure}[2][1=]{\todo[linecolor=blue,backgroundcolor=blue!25,bordercolor=blue,#1]{#2}}
\newcommandx{\changeThis}[2][1=]{\todo[linecolor=red,backgroundcolor=red!25,bordercolor=red,#1]{#2}}
\title{Stability of Discontinuous Galerkin Spectral Element Schemes for Wave Propagation when the Coefficient Matrices have Jumps}
\titlerunning{Stability of DGSEM when the Coefficient Matrices have Jumps}
\author{David A. Kopriva, Gregor J. Gassner, Jan Nordstr\"om}
\institute{David A. Kopriva \at Department of Mathematics, The Florida State University, Tallahassee, FL 32306, USA (\email{kopriva@math.fsu.edu}) and Computational Science Research Center, San Diego State University, San Diego, CA, USA\\
Gregor J. Gassner \at Department for Mathematics and Computer Science; Center for Data and Simulation Science,
University of Cologne, Weyertal 86-90, 50931, Cologne, Germany\\ 
Jan  Nordstr\"om \at Department of Mathematics, Computational Mathematics, Linköping University, 581 83 Linköping, Sweden and Department of Mathematics and Applied Mathematics
University of Johannesburg
P.O. Box 524, Auckland Park 2006, South Africa
           }
\date{}                                           

\begin{document}
\maketitle
\begin{abstract}
We use the behavior of the $L_{2}$ norm of the solutions of linear hyperbolic equations with discontinuous coefficient matrices as a surrogate to infer stability of discontinuous Galerkin spectral element methods (DGSEM). Although the $L_{2}$ norm is not bounded by the initial data for homogeneous and dissipative boundary conditions for such systems, the $L_{2}$ norm is easier to work with than a norm that discounts growth due to the discontinuities. 
We show that the DGSEM with an upwind numerical flux that satisfies the Rankine-Hugoniot (or conservation) condition has the same energy bound as the partial differential equation does in the $L_{2}$ norm, plus an added dissipation that depends on how much the approximate solution fails to satisfy the Rankine-Hugoniot jump. 
\end{abstract}

\section{Introduction}
In wave propagation problems, it is natural to find interfaces where material properties like the wave propagation speeds or density abruptly change. Examples include interfaces between two dielectrics in electromagnetic wave propagation problems, or different rock layers in geophysics. At such interfaces the solutions can make discontinuous jumps, causing difficulties for many numerical methods.

One of the key features of discontinuous Galerkin (DG) methods is that the discontinuous approximation at element interfaces naturally allows jump discontinuities in the solution if element boundaries are placed along them. Consequently, DG spectral element methods have been used for over twenty years to solve problems with material discontinuities, both stationary \cite{Koprivaetal1999},\cite{ISI:000226090600009},\cite{Hesthaven:2002uq},\cite{wilcox2010} and moving \cite{Winters:2013nx}. Computations and theory in such works show that placing the discontinuities at element boundaries leads to exponentially convergent approximations.

In a paper on discontinuous interface problems, La Cognata and Nordstr\"om \cite{La-Cognata:2016ng} noted that hyperbolic problems with discontinuous coefficients do not necessarily have their energy bounded by the initial data when measured in the $L_{2}$ norm, even with homogeneous and dissipative boundary conditions. Instead, the $L_{2}$ norm can increase or decrease, depending on the relative size of the wave speeds on either side of the discontinuity. The lack of a bound on the $L_{2}$ norm is not due to an instability in the usual sense, but is due to the fact that conservation at the interface, and the resulting jump in the solution, can increase the norm of the solution as a wave propagates across it. In an alternate norm, however, one that discounts the effect of the jump, the energy is bounded.

Here we propose a procedure where we use the $L_{2}$ norm as a surrogate to infer stability of discontinuous Galerkin spectral element methods (DGSEM) for the approximation of hyperbolic equations with discontinuous coefficient matrices. The $L_{2}$ norm is easier to work with since it does not require finding the discount factors, which are difficult to compute in general configurations of elements and interfaces. 
We show that the DGSEM with an upwind numerical flux that satisfies the Rankine-Hugoniot (or conservation) condition behaves as the partial differential equation (PDE) does in the $L_{2}$ norm, plus an added dissipation that depends on how much the approximate solution fails to satisfy the Rankine-Hugoniot jump. 

\section{Linear Hyperbolic Systems with Discontinuous Coefficients}\label{sec:ContinuousAnalysis}

In this paper we establish the stability of a discontinuous Galerkin spectral element approximation
to linear hyperbolic systems of equations of the form
\begin{equation}
\statevec u_{t} + \spacevec\nabla_x\cdot \blockvec{f} = 0,
\label{eq:DivergenceEquation}
\end{equation}
where $\statevec u$ is the state vector, and $\blockvec f$ is the vector of fluxes,
\begin{equation}
\blockvec f = \sum_{j=1}^{3}\mmatrix A_{j}\statevec u\hat x_{j} = \spacevec{\mmatrix A}\statevec u,
\end{equation}
with coefficient matrices $\mmatrix A_{j}$ and unit coordinate vectors $\hat x_{j}$. We assume throughout this paper that the coefficient matrices are piecewise
constant, with discontinuities marking what we will refer to in this paper as \textit{material interfaces}. 

\begin{figure}[htbp] 
   \centering
   \includegraphics[width=2in]{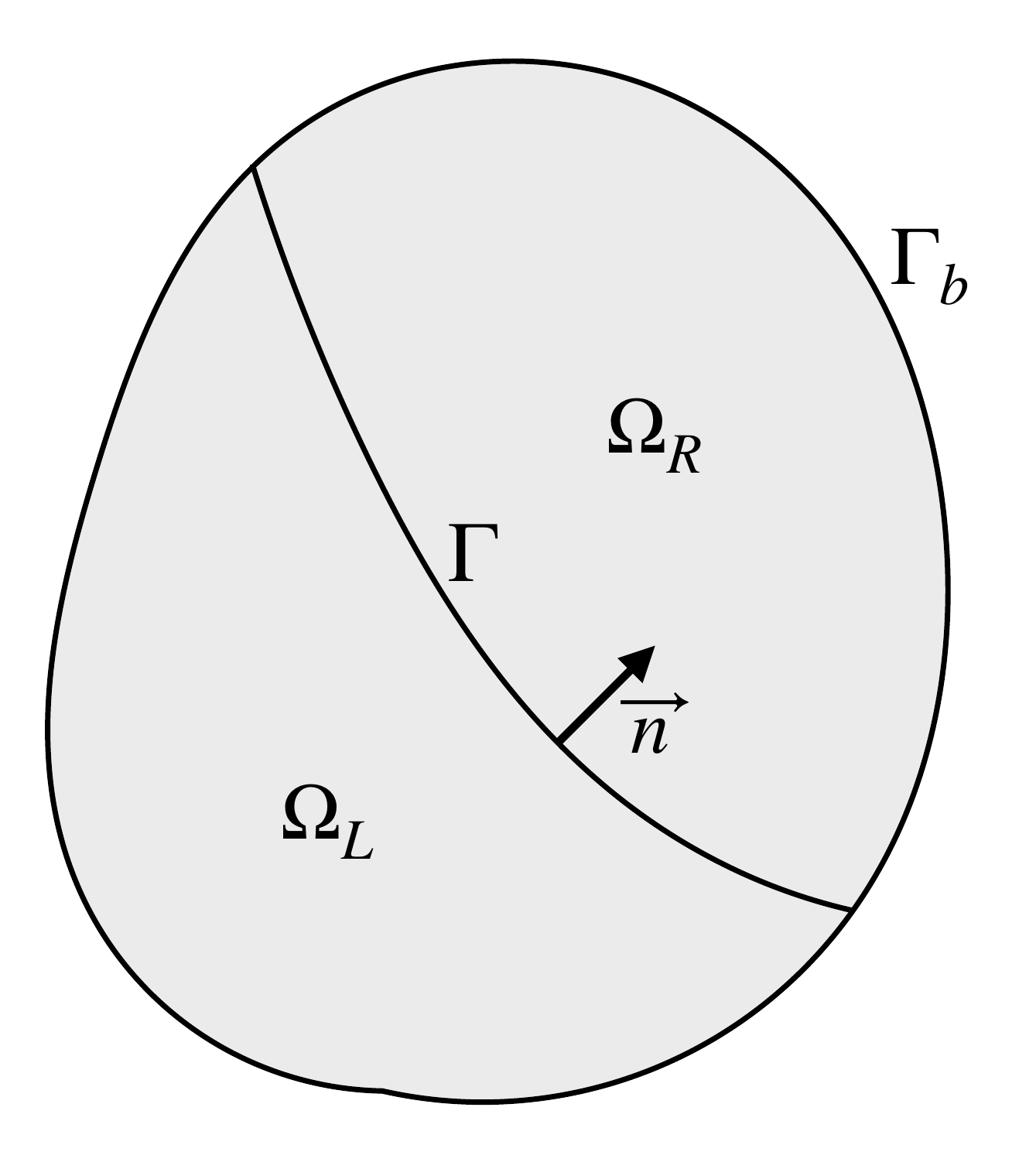} 
   \caption{Diagram of a domain $\Omega$ split by a material interface, $\Gamma$.}
   \label{fig:TwoDomains}
\end{figure}

We examine the problem defined in a domain $\Omega$, as sketched in two space dimensions in Fig. \ref{fig:TwoDomains}. It is sufficient to consider two domains with a single material interface, so the domain is split into two subdomains $\Omega_{L}$ and $\Omega_{R}$ separated by an interface $\Gamma$. The external boundary is $\Gamma_{b}$, along which we assume that proper, well-posed and dissipative boundary conditions are applied. 

Since the system is hyperbolic, there exists a matrix of right eigenvectors, $\mmatrix P$, and a real diagonal matrix, $\Lambda$, such that $\mmatrix A \equiv \spacevec k\cdot \spacevec{\mmatrix A} = \mmatrix P\Lambda \mmatrix P^{-1}$ for any nonzero space vector $\spacevec k = k_{x}\hat x + k_{y}\hat y + k_{z}\hat z$, where $(\hat x, \hat y,\hat z) = (\hat x_{1}, \hat x_{2}, \hat x_{3})$. We also assume that the matrices $\mmatrix A_{j}$
are simultaneously symmetrizable and that there exists a piecewise constant matrix $\mmatrix S$ such that 
$\mmatrix A_{j}^{s}=\mmatrix S^{-1} \mmatrix A_{j} \mmatrix S =  \left(\mmatrix A_{j}^{s}\right)^{T}$.

As a concrete example of the system \eqref{eq:DivergenceEquation}, we pose the linear acoustic wave system where
\begin{equation}
\statevec u = \left[\begin{array}{c}p \\u \\v \\w\end{array}\right],
\quad
\mmatrix A_{j} =\left[\begin{array}{cccc}0 & \delta_{j1}\rho c^2 & \delta_{j2}\rho c^2 & \delta_{j3}\rho c^2 \\\delta_{j1}/\rho & 0 & 0 & 0 \\\delta_{j2}/\rho & 0 & 0 & 0 \\\delta_{j3}/\rho & 0 & 0 & 0\end{array}\right], \quad j= 1,2,3,
\label{eq:WaveEqnMatrices}
\end{equation}
and where $\rho$ is the density of the medium, $c$ is the sound speed, and $\delta_{ij}$ is the Kronecker delta. The state vector can be viewed as representing pressure, $p$, and three velocity components, $u,v,w$.
The coefficient matrices are simultaneously symmetrizable by the matrix 
\begin{equation}
\mmatrix S = \left[\begin{array}{cccc}c & 0 & 0 & 0 \\0 & 1/\rho & 0 & 0 \\0 & 0 & 1/\rho & 0 \\0 & 0 & 0 & 1/\rho\end{array}\right].
\end{equation}
With jump discontinuities in the material parameters, $\rho$ and $c$, the coefficient matrices and symmetrizer have jump discontinuities.

We contrast the approximation of the system \eqref{eq:DivergenceEquation} with that of the approximation of systems that can be written in the form
\begin{equation}
\mmatrix E\,\statevec u_{t} + \nabla\cdot\left(\spacevec{\mmatrix B}\statevec u\right) = 0,
\label{eq:AlternatePDE}
\end{equation}
where $E>0$ is diagonal and discontinuous at material interfaces while $\spacevec{\mmatrix B}$ is continuous. The system \eqref{eq:WaveEqnMatrices}, for example,
can be re-written in the form \eqref{eq:AlternatePDE} with symmetric matrices
\begin{equation}
\mmatrix E = \left[\begin{array}{cccc}1/\rho c^{2}& 0 & 0 & 0 \\0 & \rho & 0 & 0 \\0 & 0 & \rho & 0 \\0 & 0 & 0 & \rho\end{array}\right],\quad 
\mmatrix B_{j} =\left[\begin{array}{cccc}0 & \delta_{j1} & \delta_{j2} & \delta_{j3} \\\delta_{j1} & 0 & 0 & 0 \\\delta_{j2} & 0 & 0 & 0 \\\delta_{j3} & 0 & 0 & 0\end{array}\right], \quad j= 1,2,3
\end{equation}

For equations of the form \eqref{eq:AlternatePDE}, there is a natural norm,
\begin{equation}
\inorm{\statevec u}^{2}_{E} = \int_{\Omega}\statevec u^{T}\mmatrix E\statevec u d\spacevec x,
\end{equation}
in which the energy is bounded for homogeneous dissipative physical boundary conditions and nonconservative interface conditions, with that energy satisfying
\begin{equation}
\frac{d}{dt}\inorm{\statevec u}^{2}_{E} \le 0.
\end{equation}
Stability of DG spectral approximations to equations in the form \eqref{eq:AlternatePDE} has been shown specifically, for instance, for Maxwell's equations \cite{Hesthaven:2002uq} and the elastic wave equations \cite{wilcox2010}.
\begin{remark}\label{rem:MatrixSplitRem}
The system \eqref{eq:DivergenceEquation} cannot in general be rewritten in the form \eqref{eq:AlternatePDE}. That would require that each $\mmatrix A_{j}$ can 
be written as $\mmatrix A_{j}= \mmatrix E^{-1}\mmatrix B_{j}$ where $\mmatrix E = \mmatrix E^{T}>0$ and $\mmatrix E$ contains all material properties. A counter example is the frozen coefficient compressible Euler equations \cite{Isi:A1981Lw20700001}.
\end{remark}

\begin{figure}[htbp] 
   \centering
   \includegraphics[width=4in]{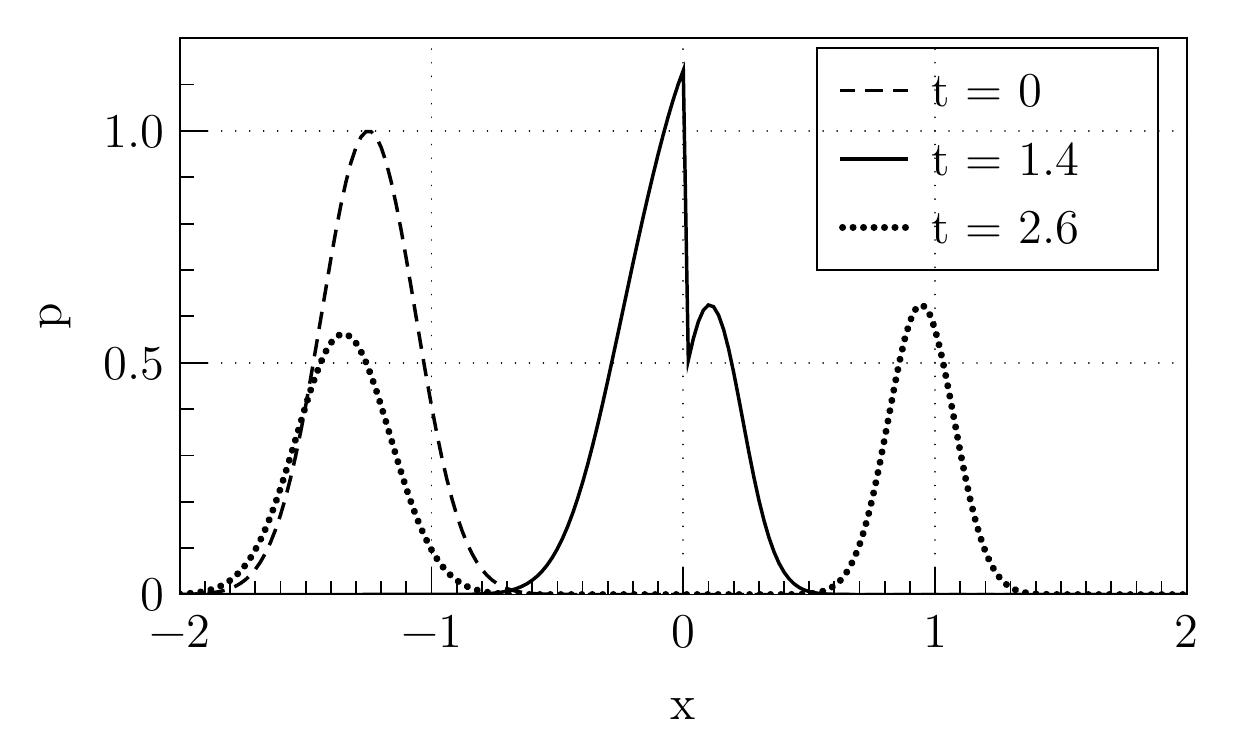} 
   \caption{Exact, analytic $p$ solution of the one dimensional acoustic wave equation for propagation of a wave across a material interface at $x=0$. The solution is plotted at three times showing the initial incident wave ($t=0$), the interaction with the material discontinuity ($t=1.4$), and the reflected and transmitted waves after the interaction ($t=2.6$)}
   \label{fig:OneDExactWaveScatter}
\end{figure}

As noted in \cite{La-Cognata:2016ng}, systems of the form \eqref{eq:DivergenceEquation} with discontinuous coefficient matrices do not necessarily have energy bounded by the initial data when measured in the $L_{2}$ norm, and we present an example here to motivate the situation. Fig. \ref{fig:OneDExactWaveScatter} shows the $p$ component of the analytic solution of acoustic wave reflection and transmission at a material boundary placed at $x = 0$ at three times: The initial incident wave, when the wave is interacting with the material discontinuity, and the reflected and transmitted waves after the interaction.

\begin{figure}[htbp] 
   \centering
   \includegraphics[width=4in]{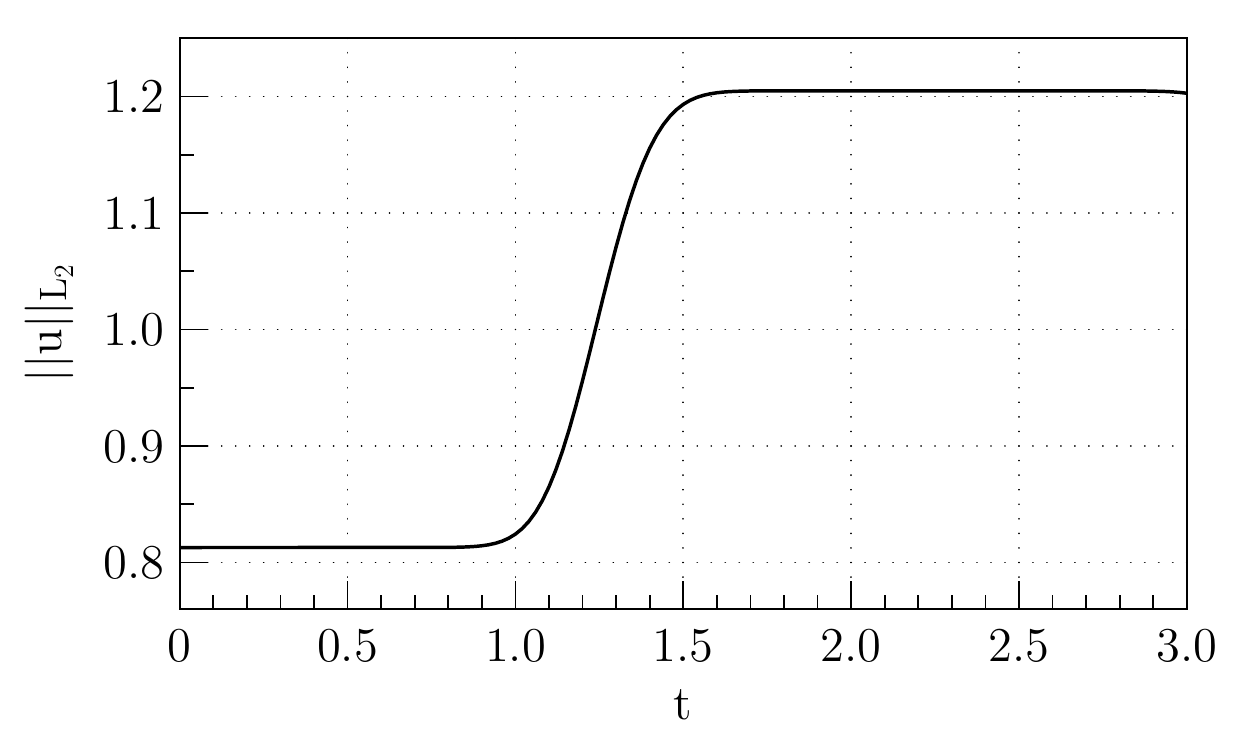} 
   \caption{Exact $L_{2}$ energy for the solution of the one dimensional acoustic wave equation for propagation of a wave across a material interface}
   \label{fig:OneDExactL2}
\end{figure}

We plot the energy as a function of time, measured by the $L_{2}$ norm,
\begin{equation}
\inorm{\statevec u}^{2}_{L_{2}} = \int_{-2}^{2}\statevec u^{T}\statevec u d x,
\end{equation}
in Fig. \ref{fig:OneDExactL2}. We see that the $L_{2}$ energy is bounded, and even though the $L_{2}$ energy estimate does not show boundedness directly, energy is bounded by the initial data in a norm that discounts the jump \cite{La-Cognata:2016ng}. Note that there is a slight downturn in the energy in Fig. \ref{fig:OneDExactL2} as $t\rightarrow 3$. The energy does decrease to zero after that time as the waves propagate out of the domain.

To establish the stability of the discontinuous Galerkin spectral element approximation of \eqref{eq:DivergenceEquation}, we follow the roadmap presented in \cite{Nordstrom:2016jk}. We first establish energy behavior of the PDE system, and then follow an equivalent discrete path to establish an equivalent behavior for the approximation. We begin with the study of the scalar one-dimensional advection problem, since it is easy to follow the steps, and then a symmetric system in one space dimension. Finally we use the symmetric system results to derive the energy bound for the general system in Sec. \ref{sec:ExtensionToDomain}.

\subsection{Energy Dynamics of the Scalar Problem in One Space Dimension}
To motivate (and simplify) the general formulation, we start with the scalar advection equation with two domains as an introduction. Our discussion in this section restates that of \cite{La-Cognata:2016ng}, but introduces our notation used in succeeding sections.

We derive the energy dynamics of the solution to the scalar advection initial-boundary-value problem in the form \eqref{eq:DivergenceEquation}
\begin{equation}
\begin{gathered}
u_{t} + au_{x}= 0 \quad x\in [-1,1]\hfill\\
u(-1,t) = 0 \hfill\\
u(x,0) = u_{0}(x),\hfill
\end{gathered}
\label{eq:ScalarIBVP}
\end{equation}
where 
\begin{equation}
a(x) = 
\left\{
\begin{gathered}
a_{L} > 0\quad x \le 0\hfill\\
a_{R} > 0\quad x > 0,\hfill
\end{gathered}
\right. 
\end{equation}
$a_{L}, a_{R}$ are constants, and $a_{L}\ne a_{R}$. The discussion that follows leads to the same types of conclusions if the wave speeds are both negative. We are interested here in problems where the domains couple and waves propagate from one side to the other. So we do not consider $a_{L}>0,a_{R}<0$, where the domains decouple as energy is dissipated at the interface, or $a_{L}<0,a_{R}>0$ where boundary conditions for both sides are required.

We split the problem into two: Left,
\begin{equation}
\begin{gathered}
u_{t }+ a_{L}u_{x} = 0 \quad x \le 0\hfill\\
u(-1,t) = 0,\hfill\\
\end{gathered}
\end{equation}
and right
\begin{equation}
\begin{gathered}
u_{t }+ a_{R}u_{x} = 0 \quad x > 0\hfill\\
u(0^{+},t) = u_{*}(t),\hfill\\
\end{gathered}
\end{equation}
where $u_{*}$ is the upwind specified interface condition chosen so that the Rankine-Hugoniot (or conservation) condition 
\begin{equation}
a_{L}u(0^{-},t) = a_{R}u_{*}(t)
\label{eq:RHScalar}
\end{equation}
is satisfied. Thus, for the scalar equation, $u_{*}(t) =\frac{a_{L}}{a_{R}} u(0^{-},t)$.

To find the energy equation, we multiply by the solution and integrate over the domains. Define the $L_{2}$ energy norms
\begin{equation}
\inorm{u}_{L}^{2}= \int_{-1}^{0} u^{2}dx,\quad \inorm{u}_{R}^{2}= \int_{0}^{1} u^{2}dx.
\end{equation}
Then
\begin{equation}
\begin{gathered}
\oneHalf\frac{d}{dt}\inorm{u}^{2}_{L} + \frac{a_{L}}{2}\left. u^{2}\right|_{-1}^{0^{-}} = 0\hfill\\
\oneHalf\frac{d}{dt}\inorm{u}^{2}_{R} + \frac{a_{R}}{2}\left. u^{2}\right|_{0^{+}}^{1} = 0.\hfill
\end{gathered}
\label{eq:BothSidesNorms}
\end{equation}
Adding together and re-arranging,
\begin{equation}
\oneHalf\frac{d}{dt}\inorm{u}^{2} -  \frac{1}{2}a_{L}u^{2}(-1) + \frac{1}{2}\left\{a_{L}u^{2}(0^{-}) -  a_{R}u^{2}(0^{+})\right\} +\frac{1}{2} a_{R}u^{2}(1)=0,
\end{equation}
where $\inorm{\cdot}^{2} = \inorm{\cdot}_{L}^{2} + \inorm{\cdot}_{R}^{2}$.
Applying the homogeneous boundary condition on the left,
\begin{equation}
\oneHalf\frac{d}{dt}\inorm{u}^{2}  + \frac{1}{2}\left\{a_{L}u^{2}(0^{-}) -  a_{R}u^{2}(0^{+})\right\}  = -\frac{1}{2} a_{R}u^{2}(1)\le 0.
\label{eq:FirstEnergyStatement}
\end{equation}
When we apply the interface condition,
\begin{equation}
\oneHalf\frac{d}{dt}\inorm{u}^{2}  \le  -\frac{1}{2}\left\{a_{L}u^{2}(0^{-}) -  a_{R}u^{2}_{*}\right\} \equiv Q.
\label{eq:ScalarEnergyDynamicsWith*}
\end{equation}
The quantity $Q$ will be used later in this paper to define stability.

Finally, we substitute the interface value for $u_{*}$,
\begin{equation}
\oneHalf\frac{d}{dt}\inorm{u}^{2}  \le  -\frac{1}{2}\left\{a_{L}u^{2}(0^{-}) -  a_{R}\frac{a_{L}^{2}}{a_{R}^{2}}u^{2}(0^{-})\right\} ,
\end{equation}
and rearrange so that
\begin{equation}
\oneHalf\frac{d}{dt}\inorm{u}^{2}  \le -\frac{a_{L}}{2}\left\{1 -  \frac{a_{L}}{a_{R}}\right\} u^{2}(0^{-}).
\label{eq:ScalarEnergyDynamics}
\end{equation}
Equation \eqref{eq:ScalarEnergyDynamics} shows that the energy is dissipated by the interface only if $a_{R} > a_{L}$. Otherwise the interface generates energy, as illustrated in Fig. \ref{fig:OneDExactL2}.

In \cite{La-Cognata:2016ng} it was shown that one can construct ``discounted norms'', in which the energy is bounded. If the second equation in \eqref{eq:BothSidesNorms} is multiplied by a constant $\alpha_{c}>0$, then the weighted sum leads to
\begin{equation}
\oneHalf\frac{d}{dt}\left\{ \inorm{u}^{2}_{L} + \alpha_{c} \inorm{u}^{2}_{R}\right\} \le -\frac{a_{L}}{2}\left\{1 -  \alpha_{c}\frac{a_{L}}{a_{R}}\right\} u^{2}(0^{-}).
\label{eq:WeightedNormTimeEnergy}
\end{equation}
Then defining the new norm with the $\alpha_{c}$ discount factor, we have
\begin{equation}
\frac{d}{dt}\inorm{u}_{\alpha_{c}}^{2}\le 0,
\label{eq:DiscountedScalarBound}
\end{equation}
provided that
\begin{equation}
\alpha_{c}\le \frac{a_{R}}{a_{L}}.
\end{equation}
\begin{remark}
The weighted norm discounts the effect of the jump, with the result that viewed in the discounted norm, the energy no longer appears to increase.
\qed\end{remark}

\begin{remark}
The use of the discounted norm scales to multiple material interfaces and multiple space dimensions by choosing $\alpha_{c}$ to be the minimum over all the ratios of downwind to upwind wave speed ratios.
\label{rem:MultipleDiscountedNorms}
\qed\end{remark}

\begin{remark}\label{rem:weightednormrem}
Alternatively, unlike the general case noted in Remark \ref{rem:MatrixSplitRem}, the scalar equation \eqref{eq:ScalarIBVP} can be recast to the form \eqref{eq:AlternatePDE} by dividing by the wave speed. Let $\varepsilon = 1/a>0$. Then 
\begin{equation}
\varepsilon u_{t}+ u_{x} = 0.
\label{eq:AltScalarEquation}
\end{equation}
If the \textit{nonconservative} boundary condition at the interface, $u_{*}=u(0^{+},t) = u(0^{-},t)$, is used, then 
following the same procedure as \eqref{eq:BothSidesNorms}--\eqref{eq:ScalarEnergyDynamics}, 
\begin{equation}
\oneHalf \frac{d}{dt}\inorm{u}_{\varepsilon}^{2}\le 0,
\end{equation}
where the weighted energy norm is given by
\begin{equation}
\inorm{u}_{\varepsilon}^{2} = \int_{-1}^{0}\varepsilon_{L}u^{2}dx + \int_{0}^{1}\varepsilon_{R}u^{2}dx.
\label{eq:AltScalarEnergy}
\end{equation}
Using the norm \eqref{eq:AltScalarEnergy}, but with the conservative interface condition \eqref{eq:RHScalar}, the solution still has a bound like \eqref{eq:ScalarEnergyDynamics}, namely
\begin{equation}
\oneHalf \frac{d}{dt}\inorm{u}_{\varepsilon}^{2}\le-\oneHalf \left\{ 1 - \left(\frac{a_{L}}{a_{R}}\right)^{2}\right\}u^{2}(0^{-}).
\end{equation} 
So when the conservative interface condition is used, the weighted energy norm is also bounded only when $a_{L}/a_{R}\le 1$ .
\qed\end{remark}

\subsection{Energy Dynamics for Hyperbolic Systems in One Space Dimension}\label{Sec:Symmetric12DsystemPDE}

We now increase the complexity and extend the scalar one-dimensional analysis to the general system \eqref{eq:DivergenceEquation} in one space dimension. We derive the energy equation for the one-dimensional hyperbolic system
\begin{equation}
\begin{split}
&\statevec u_{t} + \mmatrix A_{L}\statevec u_{x} \quad x\le 0
\\&
\statevec u_{t} + \mmatrix A_{R}\statevec u_{x} \quad x> 0,
\end{split}
\label{eq:OriginalSplit1DSystem}
\end{equation}
where the coefficient matrices are for now assumed to be symmetric. Under this assumption, there is a matrix $\mmatrix P$ such that
$ \mmatrix A =  \mmatrix P \Lambda \mmatrix P^{-1}$ 
satisfying $\mmatrix P^{-1}= \mmatrix P^{T}$. For the moment, let us assume that $\mmatrix A$ has no zero eigenvalues. We also assume that the number of positive and negative eigenvalues does not change across the interface. In other words, there is no eigenvalue that changes sign at the jump. Depending on the sign change, boundary/interface conditions are either lost or gained. More general conditions where the sign of the eigenvalues changes in multi-physics applications are considered in \cite{doi:10.1137/16M1087710}. Finally, we assume that appropriate boundary and initial data are applied.

To find the interface condition at $x=0$ for the system \eqref{eq:OriginalSplit1DSystem}, we split the system into right and left going waves. The characteristic variables for the system \eqref{eq:OriginalSplit1DSystem} are 
\begin{equation} 
\statevec w = \mmatrix P^{-1}\statevec u = \left[\begin{array}{c}\statevec w^{+} \\\statevec w^{-}\end{array}\right],
\end{equation}
where $\statevec w^{+}$ is associated with the positive eigenvalues of $\mmatrix A$ and $\statevec w^{-}$ is associated with the negative ones.
They are chosen upwind at the interface according to
\begin{equation}
\statevec w^{+}_{R}=\statevec w^{+}_{*},\quad \statevec w^{-}_{L}=\statevec w^{-}_{*},
\label{eq:UpwindedCharVars}
\end{equation}
where here and in the following, the subscripts $R$ and $L$ correspond to the values at $x=0^{+}$ and $x=0^{-}$, respectively.

The $\statevec w^{\pm}_{*}$ are computed so that the Rankine-Hugoniot condition
 \begin{equation}
 \mmatrix A_{L}\statevec u|_{0^{-}} = \mmatrix A_{R}\statevec u|_{0^{+}}\Leftrightarrow
 \mmatrix P_{L} \Lambda_{L}\left[\begin{array}{c}\statevec w_{L}^{+} \\\statevec w_{*}^{-}\end{array}\right] =\mmatrix P_{R} \Lambda_{R}\left[\begin{array}{c}\statevec w_{*}^{+} \\\statevec w_{R}^{-}\end{array}\right]
 \label{eq:RHCondition1DSystem}
 \end{equation}
  is satisfied at the stationary interface. 
Let us write
 \begin{equation}
 \Lambda = \left[\begin{array}{cc}\bar \Lambda^+ & 0 \\0 & \bar \Lambda^-\end{array}\right],\quad
  \Lambda^{+} = \left[\begin{array}{cc}\bar \Lambda^+ & 0 \\0 & 0\end{array}\right],\quad \Lambda^{-} = \left[\begin{array}{cc}0 & 0 \\0 & \bar \Lambda^-\end{array}\right].
 \end{equation} 
 Then \eqref{eq:RHCondition1DSystem} can be written as
 \begin{equation}
\mmatrix  P_{L} \Lambda^{+}_{L}\left[\begin{array}{c}\statevec w_{L}^{+} \\0\end{array}\right] 
 + \mmatrix P_{L} \Lambda^{-}_{L}\left[\begin{array}{c}\statevec 0\\\statevec w_{*}^{-}\end{array}\right] 
  =
   \mmatrix P_{R} \Lambda^{+}_{R}\left[\begin{array}{c}\statevec w_{*}^{+} \\0\end{array}\right] 
 + \mmatrix P_{R} \Lambda^{-}_{R}\left[\begin{array}{c}\statevec 0\\\statevec w_{R}^{-}\end{array}\right] .
 \end{equation}
 Let us put the unknowns on the left, and the knowns on the right, giving
 \begin{equation}
 \mmatrix P_{L} \Lambda^{-}_{L}\left[\begin{array}{c}\statevec 0\\\statevec w_{*}^{-}\end{array}\right]  - 
   \mmatrix  P_{R} \Lambda^{+}_{R}\left[\begin{array}{c}\statevec w_{*}^{+} \\0\end{array}\right] =
   \mmatrix  P_{R} \Lambda^{-}_{R}\left[\begin{array}{c}\statevec 0\\\statevec w_{R}^{-}\end{array}\right] -
\mmatrix  P_{L} \Lambda^{+}_{L}\left[\begin{array}{c}\statevec w_{L}^{+} \\0\end{array}\right] .
 \label{eq:wstareqn1}
  \end{equation}
 Eq. \eqref{eq:wstareqn1} provides a system of equations for the unknowns.
 
 The matrices on the left of \eqref{eq:wstareqn1} have a special structure since $\mmatrix P$ is the matrix of right eigenvectors and $\Lambda$ is a diagonal matrix. Let $n$ be the number of positive eigenvalues out of a total of $m$.
 Then define
\begin{equation}
\mmatrix M^{+}\equiv \mmatrix P\Lambda^{+}=\left[\begin{array}{cccccc}\lambda_1 \spacevec p_1 & \ldots & \lambda_n \spacevec p_n & 0 &\ldots& 0\end{array}\right]
\end{equation}
and
\begin{equation}
\mmatrix M^{-}\equiv \mmatrix P\Lambda^{-} = \left[\begin{array}{cccccc}0 & \ldots & 0 & \lambda_n \spacevec p_{n+1} & \ldots &  \lambda_{m}\spacevec p_{m}\end{array}\right],
\end{equation}
where  $\spacevec p_{j}$ is the eigenvector associated with the eigenvalue $\lambda_{j}$ and the eigenvalues are ordered in decreasing order, largest to smallest with $\lambda_{j}>0$ for $j\le n$. Then we can write \eqref{eq:wstareqn1} as
 \begin{equation}
 \mmatrix M_{L}^{-}\left[\begin{array}{c}\statevec 0\\\statevec w_{*}^{-}\end{array}\right]  - 
   \mmatrix M_{R}^{+}\left[\begin{array}{c}\statevec w_{*}^{+} \\0\end{array}\right] =
    \mmatrix M_{R}^{-}\left[\begin{array}{c}\statevec 0\\\statevec w_{R}^{-}\end{array}\right] -
 \mmatrix M_{L}^{+}\left[\begin{array}{c}\statevec w_{L}^{+} \\0\end{array}\right] .
 \label{eq:wstareqn2}
  \end{equation}
 Given the structure of the $\mmatrix M^{\pm}$ matrices, the equations can be combined to produce a single system for the unknowns
  \begin{equation}
 \mmatrix M_{LR}\left[\begin{array}{c}\statevec w_{*}^{+}\\\statevec w_{*}^{-}\end{array}\right]  
  =
    \mmatrix M_{RL}\left[\begin{array}{c}\statevec w_{L}^{+}\\\statevec w_{R}^{-}\end{array}\right],
 \label{eq:wstareqn3}
  \end{equation}
where 
\begin{equation}
\mmatrix M_{LR}\equiv \mmatrix M^{-}_{L} - \mmatrix M^{+}_{R},\quad \mmatrix M_{RL}\equiv \mmatrix M^{-}_{R} - \mmatrix M^{+}_{L}.
\end{equation}

Existence and uniqueness of the inflow characteristic vectors $\statevec w_{*}^{\pm}$ therefore depends on the existence of the inverse of the matrix $\mmatrix M_{LR}$. That matrix is comprised of eigenvectors of the coefficient matrix evaluated on the left and eigenvectors evaluated on the right.
On the one hand, if the eigenvectors of the coefficient matrix do not change across the material discontinuity, then, since the eigenvectors are independent, $\mmatrix M_{LR}^{-1} \mmatrix M_{LR}$ is diagonal. As an example, the eigenvectors of the acoustic wave system \eqref{eq:WaveEqnMatrices} are constant, being independent of the material properties on either side. On the other hand, if the eigenvectors change across the interface and the matrix $\mmatrix M_{LR}^{-1}\mmatrix M_{RL}$ is not diagonal, then the problem is ill-posed \cite{doi:10.1137/16M1087710}. 
We therefore require that the eigenvectors be preserved across the jumps so that $\mmatrix M_{LR}^{-1}$ exists,  $\mmatrix M \equiv \mmatrix M_{LR}^{-1} \mmatrix M_{LR}$ is diagonal, and 
  \begin{equation}
 \left[\begin{array}{c}\statevec w_{*}^{+}\\\statevec w_{*}^{-}\end{array}\right]  
  =
    \mmatrix M_{LR}^{-1}\mmatrix M_{RL}\left[\begin{array}{c}\statevec w_{L}^{+}\\\statevec w_{R}^{-}\end{array}\right] \equiv \mmatrix M\left[\begin{array}{c}\statevec w_{L}^{+}\\\statevec w_{R}^{-}\end{array}\right].
 \label{eq:wstarSoln}
  \end{equation}

\begin{remark}
The Rankine-Hugoniot (conservation) condition \eqref{eq:RHCondition1DSystem} limits the form of the interface condition significantly. If only boundedness is desired, more general coupling conditions are allowed \cite{doi:10.1137/16M1087710}.
\qed\end{remark}
\begin{remark}
 One can see that the assumption that there are no zero eigenvalues is not a restriction. If there are zero eigenvalues, then the associated characteristic variables $\statevec w^{0}$ are multiplied by the zero matrix and have no contribution to the system. Therefore those quantities can be eliminated, leaving \eqref{eq:wstareqn2} and what followed. The $\statevec w^{0}$ vector is determined by the initial data.
 \qed\end{remark}

Going back to the original equations, \eqref{eq:OriginalSplit1DSystem}, we compute the energy equation by multiplying by the state and integrating over the domain, giving
\begin{equation}
\oneHalf\frac{d}{dt} \left\{ \inorm{\statevec u}_{L}^{2}+  \inorm{\statevec u}_{R}^{2}\right\} + \PBT =-\oneHalf \left\{\statevec u^{T}_{L}  \mmatrix A_{L}\statevec u_{L} - \statevec u^{T}_{R} \mmatrix A_{R}\statevec u_{R}\right\} ,
\label{eq:System1DEnergy}
\end{equation}
where, now, $\inorm{\statevec u}^{2} = \iprod{\statevec u,\statevec u}$ and $\PBT$ represents the terms coming from the physical boundary conditions on the left and right. Since we are only interested here in the interface conditions, we will assume that the physical boundary conditions are well posed so that $\PBT \ge 0$. In that case, 
\begin{equation}
\oneHalf\frac{d}{dt} \left\{ \inorm{\statevec u}_{L}^{2}+  \inorm{\statevec u}_{R}^{2}\right\} \le Q ,
\label{eq:OneDSystemNormBound}
\end{equation}
where
\begin{equation}
Q \equiv -\oneHalf\left\{\statevec u^{T}_{L} \mmatrix A_{L}\statevec u_{L} - \statevec u^{T}_{R}  \mmatrix A_{R}\statevec u_{R}\right\}
\end{equation}
is the interface contribution to the energy.

Following the steps in the scalar analysis, we now apply the interface boundary conditions on $Q$. We decompose the system into characteristic variables. 
Then we use the fact that $\mmatrix A$ is symmetric, making  $\mmatrix P^{-1}= \mmatrix P^{T}$. With this decomposition,
\begin{equation}
\begin{split}
Q &= -\oneHalf\left\{\statevec w_{L}^{T} \Lambda_{L}\statevec w_{L} - w_{R}^{T}\Lambda_{R} \statevec w_{R}\right\}
\\&
= -\oneHalf\left\{\statevec w^{+,T}_{L} \bar\Lambda^{+}_{L}\statevec w^{+}_{L} + \statevec w^{-,T}_{*} \bar\Lambda^{-}_{L}\statevec w^{-}_{*}\right\}
 +\oneHalf \left\{\statevec w^{+,T}_{*} \bar\Lambda^{+}_{R}\statevec w^{+}_{*} + \statevec w^{-,T}_{R} \bar\Lambda^{-}_{R}\statevec w^{-}_{R}\right\},
\end{split}
\end{equation}
taking into account the upwinding of the characteristic variables, \eqref{eq:UpwindedCharVars}.

We now gather the right-going and left-going wave contributions (c.f. \eqref{eq:ScalarEnergyDynamicsWith*}), 
\begin{equation}
\begin{split}
Q 
&=  -\oneHalf\left\{\statevec w^{+,T}_{L} \bar\Lambda^{+}_{L}\statevec w^{+}_{L} - \statevec w^{+,T}_{*} \bar\Lambda^{+}_{R}\statevec w^{+}_{*}\right\}
 + \oneHalf\left\{  \statevec w^{-,T}_{R} \bar\Lambda^{-}_{R}\statevec w^{-}_{R} - \statevec w^{-,T}_{*} \bar\Lambda^{-}_{L}\statevec w^{-}_{*}\right\},
\end{split}
\label{eq:QForSystem0}
\end{equation}
and then use the fact that $\bar\Lambda^{-}<0$, to get the final form of the interface contribution, which we write in terms of its characteristic components,
\begin{equation}
Q\left(\statevec w_{L}, \statevec w_{R}\right)
 = -\oneHalf\left\{\statevec w^{+,T}_{L} \bar\Lambda^{+}_{L}\statevec w^{+}_{L} - \statevec w^{+,T}_{*} \bar\Lambda^{+}_{R}\statevec w^{+}_{*}\right\}
 - \oneHalf\left\{  \statevec w^{-,T}_{R}\left|\bar\Lambda^{-}_{R}\right|\statevec w^{-}_{R} - \statevec w^{-,T}_{*} \left|\bar\Lambda^{-}_{L}\right|\statevec w^{-}_{*}\right\}.
\label{eq:QForSystem}
\end{equation}
Eq. \eqref{eq:QForSystem} is the system version of the scalar interface condition seen in \eqref{eq:ScalarEnergyDynamicsWith*}. 

As in the scalar problem, one can construct a discounted norm $\inorm{\cdot}_{B}$ for which the associated interface term $Q_{B}$ is non-positive and the discounted norm is bounded when the coupling matrix, $\mmatrix M$, exists and is diagonal. For instance, one simple choice is to let
\begin{equation}
\mmatrix B = \mmatrix P\left[\begin{array}{cccccc}\mu^+ &  &  &  &  &  \\ & \ddots &  &  &  &  \\ &  & \mu^+ &  &  &  \\ &  &  & \mu^- &  &  \\ &  &  &  & \ddots &  \\ &  &  &  &  & \mu^-\end{array}\right]\mmatrix P^{-1},
\end{equation}
where the entries with $\mu^{\pm}>0$ are counted according to the number of positive and negative eigenvalues of $\mmatrix A$. Then multiplying the system on $x>0$ by $\mmatrix B$ from the left, defining the norm $\inorm{\statevec u}_{B} = \iprod{\statevec u, \mmatrix B\statevec u}^{\oneHalf}$, and following the same steps leading to \eqref{eq:QForSystem}, the interface contribution to the energy is
\begin{equation}
Q_{B}
 = -\oneHalf\left\{\statevec w^{+,T}_{L} \bar\Lambda^{+}_{L}\statevec w^{+}_{L} - \mu^{+}\statevec w^{+,T}_{*} \bar\Lambda^{+}_{R}\statevec w^{+}_{*}\right\}
 - \oneHalf\left\{  \mu^{-}\statevec w^{-,T}_{R}\left|\bar\Lambda^{-}_{R}\right|\statevec w^{-}_{R} - \statevec w^{-,T}_{*} \left|\bar\Lambda^{-}_{L}\right|\statevec w^{-}_{*}\right\}
\label{eq:QSubBForSystem}.
\end{equation}

One then only needs to find $\mu^{+}$ small enough and $\mu^{-}$ large enough to ensure that $Q_{B}\le 0$, in which case the new energy $\sqrt{\inorm{\statevec u}^{2}_{L} + \inorm{\statevec u}^{2}_{B,R}}$ is bounded by the initial data. Since the coupling matrix $\mmatrix M$ is diagonal, let us split it as
\begin{equation}
\mmatrix M = \left(\begin{array}{cc}\bar M^+ & 0 \\0 & \bar M^-\end{array}\right)
\end{equation}
so that $\statevec w^{+}_{*}= \bar M^+\statevec w^{+}_{L}$ and $\statevec w^{-}_{*}= \bar M^-\statevec w^{-}_{R}$. Then $Q_{B}\le 0$ if $\mu^{\pm}$ are chosen so that
\begin{equation}
\begin{split}
\bar\Lambda_{L}^{+} - \mu^{+}\bar M^{+,T}\Lambda_{R}^{+}\bar M^{+}&>0,\\
\mu^{-}\left|\bar\Lambda_{R}^{-}\right| - \bar M^{-,T}\left|\bar\Lambda_{L}^{-}\right| \bar M^{-}&>0.
\end{split}
\label{eq:MuEquations}
\end{equation}

\subsection{Extension to Non-Symmetric Equations and an Arbitrary Domain }\label{sec:ExtensionToDomain}

The results of the previous two sections extend to general geometries and non-symmetric coefficient matrices. 
In preparation for the generalization of \eqref{eq:QForSystem}, we note that within each subdomain the coefficient matrices are constant, and therefore we can re-write \eqref{eq:DivergenceEquation} in a split form as 
\begin{equation}
\statevec u_{t} + \oneHalf\left\{ \nabla\cdot\left(\spacevec {\mmatrix A}\statevec u\right) + \spacevec {\mmatrix A}\cdot\nabla\statevec u\right\} = 0.
\label{eq:SplitFormPDE}
\end{equation}
We also define the inner product and norm over a subdomain $D = \Omega_{L}$ or $ \Omega_{R}$ as
\begin{equation}
\iprod{\statevec u, \statevec v}_{D} = \int_{D}\statevec u^{T} \statevec vd\statevec x, \quad \inorm{\statevec u}_{D} = \iprod{\statevec u, \statevec u}_{D}^{\oneHalf}
\end{equation}
so that
\begin{equation}
 \inorm{\statevec u}^{2}_{\Omega} = \inorm{\statevec u}^{2}_{\Omega_{L}} + \inorm{\statevec u}^{2}_{\Omega_{R}}.
\end{equation}

To form the energy, we take the inner product of \eqref{eq:SplitFormPDE} with the vector $\left(\mmatrix S^{-1}\right)^{T}\mmatrix S^{-1}\statevec u$, giving
\begin{equation}
\iprod{\left(\mmatrix S^{-1}\right)^{T}\mmatrix S^{-1}\statevec u,\statevec u_{t}} _{D}
+ \oneHalf\iprod{\left(\mmatrix S^{-1}\right)^{T}\mmatrix S^{-1}\statevec u, \nabla\cdot\left(\spacevec {\mmatrix A}\statevec u\right)}_{D} + \oneHalf\iprod{\left(\mmatrix S^{-1}\right)^{T}\mmatrix S^{-1}\statevec u,  \spacevec {\mmatrix A}\cdot\nabla\statevec u}_{D} = 0.
\end{equation}
Let us define $ \sym{\statevec u} = S^{-1}\statevec u$ to be the symmetric system state. Then since $\mmatrix S$ is constant within the subdomains and $\mmatrix A^{s}=\mmatrix S^{-1} \mmatrix A \mmatrix S $,
\begin{equation}
\oneHalf\frac{d}{dt}\inorm{\sym{\statevec u}}_{D}^{2} 
+\oneHalf\iprod{\sym{\statevec u}, \nabla\cdot\left(\sym{\spacevec {\mmatrix A}}\sym{\statevec u}\right)}_{D} 
+ \oneHalf\iprod{\sym{\statevec u},  \sym{\spacevec {\mmatrix A}}\cdot\nabla\sym{\statevec u}}_{D} = 0.
\label{eq:SymmEEqn1}
\end{equation}
We then apply multidimensional integration by parts and symmetry to the divergence term 
\begin{equation}
\iprod{\sym{\statevec u}, \nabla\cdot\left(\sym{\spacevec {\mmatrix A}}\sym{\statevec u}\right)}_{D} 
= \int_{\partial D}\statevec u^{s,T}{\spacevec {\mmatrix A}^{s}}\cdot \spacevec n\sym{\statevec u}\dS - \iprod{\sym{\statevec u},  \sym{\spacevec {\mmatrix A}}\cdot\nabla\sym{\statevec u}}_{D} ,
\end{equation}
where $\spacevec n$ is the outward normal at the boundary of $D$, and note that the
volume term cancels the third term in \eqref{eq:SymmEEqn1}, leaving only the boundary integral,
\begin{equation}
\oneHalf\frac{d}{dt}\inorm{\sym{\statevec u}}_{D}^{2}  = - \oneHalf\int_{\partial D}\statevec u^{s,T}{\spacevec {\mmatrix A}^{s}}\cdot \spacevec n\sym{\statevec u}\dS .
\end{equation}
Then over the domain $\Omega$,
\begin{equation}
\oneHalf\frac{d}{dt}\inorm{\sym{\statevec u}}_{\Omega}^{2} = -\oneHalf\int_{\Gamma_{b}}\statevec u^{s,T}\left({\spacevec {\mmatrix A}^{s}}\cdot \spacevec n\right)\sym{\statevec u}\dS
 - \oneHalf\int_{\Gamma}\left\{\statevec u_{L}^{s,T}{\spacevec {\mmatrix A}_{L}^{s}}\cdot \hat n\sym{\statevec u}_{L}-\statevec u_{R}^{s,T}{\spacevec {\mmatrix A}_{R}^{s}}\cdot \hat n\sym{\statevec u}_{R}\right\}\dS,
\end{equation}
where $L/R$ represent the states on either side of the interface with respect to the normal $\spacevec n$. 

We can now get a bound for the multidimensional system similar to \eqref{eq:OneDSystemNormBound}. The integrand in the interface integral is identical to that in \eqref{eq:System1DEnergy}, with $A \leftarrow {\spacevec {\mmatrix A}_{L}^{s}}\cdot \hat n$ and $\statevec u \leftarrow \sym{\statevec u}$.  Therefore, if the boundary conditions along $\Gamma_{b}$ are properly posed and dissipative, 
\begin{equation}
\oneHalf\frac{d}{dt}\inorm{\sym{\statevec u}}_{\Omega}^{2}\le \int_{\Gamma}Q\dS.
\label{eq:GeneralWell-posednessCondition}
\end{equation}
Therefore, $Q$ is still given by \eqref{eq:QForSystem}, but now formulated in the new symmetrized variables. Note that the norm defined by $\inorm{\sym{\statevec u}}^{2} = \iprod{\left(\mmatrix S^{-1}\right)^{T}\mmatrix S^{-1}\statevec u,\statevec u}$ is equivalent to the norm $\inorm{\statevec u}$ since $\left(\mmatrix S^{-1}\right)^{T}\mmatrix S^{-1}>0$. 

\subsection{Stability}

In summary, for hyperbolic systems of the form \eqref{eq:DivergenceEquation}, with discontinuities in the coefficient matrices and homogeneous, dissipative boundary conditions, the $L_{2}$ norm of the solution obeys \eqref{eq:GeneralWell-posednessCondition}. The integrand of the interface contribution, $Q$, is of the form \eqref{eq:QForSystem}, where the characteristic variables are evaluated from the upwind side and satisfy the Rankine-Hugoniot condition. It is not necessarily non-negative, depending on the relative wave speeds from either side of the interface, so the $L_{2}$ norm of the solution is not bounded in general by the initial data. An example of such behavior was shown in Fig. \ref{fig:OneDExactL2}.

Although the $L_{2}$ norm (or, for that matter, weighted norms, see Remark \ref{rem:weightednormrem}) is not always bounded by the initial data, there exists an energy in a discounted norm that is bounded in the usual way provided that the coupling matrix between the upwind and downwind states is diagonal. 

Thus, we have two views of stability at our disposal, which we will call direct and inferred: 
\begin{itemize}
\item \textbf{Direct Stability}. When the $L_{2}$ norm is bounded, we directly have $L_{2}$ stability. This is seen in scalar problems if $a_{L} / a_{R}\le 1$ in \eqref{eq:ScalarEnergyDynamics}. For the system, the equivalent is when $\bar\Lambda_{L}^{+} - \bar M^{+,T}\Lambda_{R}^{+}\bar M^{+}>0$ and 
$\left|\bar\Lambda_{R}^{-}\right| - \bar M^{-,T}\left|\bar\Lambda_{L}^{-}\right| \bar M^{-}>0$, as seen through \eqref{eq:MuEquations} setting $\mu^{\pm} = 1$.
\item \textbf{Inferred Stability}. Nonetheless, even if the $L_{2}$ norm is not directly bounded, we have seen that one can construct a discounted norm in which it is, e.g. \eqref{eq:DiscountedScalarBound} for scalar problems and for systems when \eqref{eq:MuEquations} is satisfied. Stability in some discounted norm is therefore inferred, or implicit, if $Q$ is given by \eqref{eq:QForSystem}.
\end{itemize}

In general geometries it may not be easy to find the discounted norm in which the solution is bounded. Finding the precise coefficients requires satisfying conditions like \eqref{eq:MuEquations}. When multiple subdomains exist in multiple space dimensions, $T$-type intersections between materials are possible. The discount factors must then take into account all subdomain boundaries and be adjusted globally so that at each interface \eqref{eq:MuEquations} is still satisfied. 
For these reasons, it is easier to monitor the behavior \eqref{eq:GeneralWell-posednessCondition} of the simpler $L_{2}$ norm as a surrogate to infer well-posedness of the system. Further insights into how choosing the norm affects how the energy is bounded or not can be found in \cite{Manzanero:2018kk}.

Stability of a numerical approximation of a system follows that of the PDE, and so we state the stability condition for the approximation as:
\begin{definition}
\label{def:StabilityDefinition}
A scheme approximating the discontinuous coefficient problem \eqref{eq:DivergenceEquation} is said to have inferred stability if the discrete approximation of the standard $L_{2}$ norm is bounded as in \eqref{eq:GeneralWell-posednessCondition} and the approximation to the integrand, $Q_{N} \approx Q$, satisfies
\[
Q_{N} \le Q(\statevec W_{L},\statevec W_{R}),
\]
where $\statevec W$ is the approximation of $\statevec w$.
\end{definition}

\section{The discontinuous Galerkin spectral element discretization}

In this section, we briefly summarize the important discretization steps. For a detailed description and derivation of the scheme, we refer to e.g. \cite{Kopriva:2009nx,Gassner_BR1,winters2020construction}. 

The first step is to divide the computational domain into a mesh of non-overlapping, possibly curved, hexahedral (quadrilateral in 2D) elements, $\{e^l\}_{l=1}^K$. Each hexahedron is mapped from physical space to a reference space cube $E = [-1,1]^3$ with $x = \spacevec{X}^l(\spacevec\xi)$. From the mapping, we can compute the metric terms
\begin{equation}
\spacevec{a}_i = \frac{\partial\spacevec{X}}{\partial\xi_i},\quad i=1,2,3;\quad J = \spacevec{a}_1\cdot (\spacevec{a}_2\times \spacevec{a}_3);\quad J\spacevec{a}^i= \spacevec{a}_j\times \spacevec{a}_k,\quad (i,j,k)\text{ cyclic}.
\label{eq:metric_stuff}
\end{equation} 
Note that we need to carefully evaluate the metric terms to get a discretely divergence-free contravariant basis $J\spacevec{a}^i$, which is necessary to guarantee free-stream preservation of the discretization \cite{Kopriva:2006er} and stability of the volume terms \cite{Kopriva2016274,Gassner_BR1}.

The second step of the discetization process is to transform the problem \eqref{eq:SplitFormPDE} from physical to reference space. In reference space, \eqref{eq:SplitFormPDE} becomes
\begin{equation}
J\,\statevec u_{t} + \frac{1}{2}\left\{ \spacevec{\nabla}_\xi\cdot\left(\bigmatrix{M}^T\,\spacevec {\mmatrix A}\statevec u\right) + \spacevec {\mmatrix A}\cdot\bigmatrix{M}\,\spacevec{\nabla}_\xi\statevec u\right\} = 0,
\label{eq:split_transformed}
\end{equation}
where we collect the metric terms in the block matrix 
\begin{equation}
\bigmatrix{M} = 
\begin{pmatrix}
Ja_1^1\,\mmatrix{I} & Ja_1^2\,\mmatrix{I} & Ja_1^3\,\mmatrix{I}\\
Ja_2^1\,\mmatrix{I} & Ja_2^2\,\mmatrix{I} & Ja_2^3\,\mmatrix{I}\\
Ja_3^1\,\mmatrix{I} & Ja_3^2\,\mmatrix{I} & Ja_3^3\,\mmatrix{I}\\
\end{pmatrix},
\end{equation}
 with the identity matrix, $\mmatrix{I}$, having the size as the state vector $\statevec u$.
 
The third step is the variational Galerkin formulation. We first approximate the solution with an interpolatory polynomial of degree $N$, and denote polynomial approximations with capital letters ${u\approx U = \mathbb{I}^N(u)}$, where $\mathbb{I}^N$ denotes the interpolation operator. In the spectral collocation framework, one typically uses a nodal basis for the interpolation. Furthermore, for hexahedral/quadrilateral elements, we use a tensor-product of one-dimensional nodal Lagrange basis functions  spanned on the Legendre-Gauss-Lobatto nodes. The same polynomial approximation is used for all quantities, e.g. for the contravariant flux function $\blockvec{\widetilde{f}}\approx \blockvec{\widetilde{F}} = \mathbb{I}^N(\bigmatrix{M}^T\blockvec f)$. 

To get the variational formulation, we multiply the transformed PDE \eqref{eq:split_transformed} by polynomial test functions $\statevecGreek{\varphi}$, which are linear combinations of the nodal basis functions. Then we integrate over the reference element $E$ and use integration-by-parts to arrive at 
\begin{equation}
\begin{split}
\iprod{\mathbb{I}^N(J)\,\statevec U_t,\statevecGreek\varphi}_E &- \frac{1}{2} \iprod{\mathbb{I}^N\left(\bigmatrix{M}^T\,\spacevec {\mmatrix A}\statevec U\right),\spacevec{\nabla}_\xi\,\mathbf\varphi}_E  + \frac{1}{2} \iprod{\mathbb{I}^N\left(\spacevec {\mmatrix A}\cdot\bigmatrix{M}\,\spacevec{\nabla}_\xi\statevec U\right),\statevecGreek\varphi}_E 
\\&= - \int\limits_{\partial E} \statevecGreek\varphi^T\left\{\blockvec{\widetilde{F}} -  \frac{1}{2}\mathbb{I}^N\left(\bigmatrix{M}^T\,\spacevec {\mmatrix A}\statevec U\right)\right\}\cdot\hat{n}\,\dS,
\end{split}
\label{eq:DGWithExactIntegration}
\end{equation}
where $\hat{n}$ is the reference space outward pointing normal vector to the face $\partial E$.

Finally, we replace the integration in \eqref{eq:DGWithExactIntegration} by quadrature and cubature rules, collocated with the Legendre-Gauss-Lobatto interpolation. Note, that the Legendre-Gauss-Lobatto nodes include the boundary nodes and hence surface and volume integration nodes partially coincide with the interpolation ansatz and $\mathbb{I}^N(\cdot)$ can be dropped. Furthermore, we introduce the yet to be defined numerical flux function $\statevec{F}_n^*=\statevec{F}_n^*(\statevec U^L,\statevec U^R)\approx \blockvec{\widetilde{F}}\cdot\spacevec{n}$, which depends on the two states $\statevec{U}^{L,R}$ at the interface and approximates the normal flux through the interface. Note that we assume the coefficients ${\mmatrix A}$ are mostly constant, but when they jump, the mesh is aligned so that an element interface is at the jump. Hence, the numerical flux function at the coefficient jump interface depends not only on the solutions left and right, but also on the coefficients left and right: $\statevec{F}_n^*=\statevec{F}_n^*(\statevec U^{L,R};{\mmatrix A}^{L,R})$.

Applying quadrature, we get the formal statement of the DGSEM,
\begin{equation}
\iprod{J\,\statevec U_t,\statevecGreek\varphi}_N - \frac{1}{2} \iprod{\bigmatrix{M}^T\,\spacevec {\mmatrix A}\statevec U,\spacevec{\nabla}_\xi\,\statevecGreek\varphi}_N  + \frac{1}{2} \iprod{\spacevec {\mmatrix A}\cdot\bigmatrix{M}\,\spacevec{\nabla}_\xi\statevec U,\statevecGreek\varphi}_N = - \int\limits_{\partial E,N} \statevecGreek\varphi^T\left\{\statevec{F}_n^* -  \frac{1}{2}\statevec{F}_n\right\}\,\dS,
\label{eq:DGSEM}
\end{equation}
where $\iprod{\cdot,\cdot}_{N}$ and $\int\limits_{\partial E,N}$ represent the volume and surface quadratures, see \cite{10.1007/978-3-319-65870-4_2}. The right hand side of \eqref{eq:DGSEM} is written in terms of the normal covariant fluxes and is equivalent to that written in terms of the contravariant ones \cite{winters2020construction}. 
The resulting high-order semi-discretization is integrated with a proper high-order accurate explicit Runge-Kutta time integrator, which is stable under the typical CFL-type time step restriction.
  
\section{Stability of the Discontinuous Galerkin Approximation}
We establish the stability bound from the weak form of the equation, \eqref{eq:DGSEM}. We then follow the path taken in Sec. \ref{sec:ContinuousAnalysis} for the continuous problem to examine the discontinuous interface term:  We examine the scalar problem for insights, then the symmetric one-dimensional system, and finally the general problem for the DGSEM approximation.
\subsection{Discrete stability estimate}
For a detailed derivation of the discrete stability estimate, which parallels the continuous analysis, we refer to \cite{Gassner_BR1,winters2020construction}. Here, we will only sketch  some important intermediate steps. To get the stability estimate, we replace the test function $\statevecGreek\varphi$ with the approximate solution polynomial and the symmetrizer matrices, writing $\statevecGreek\varphi = \left(\mmatrix S^{-1}\right)^{T}\mmatrix S^{-1}\statevec U = \left(\mmatrix S^{-1}\right)^{T}\statevec U^s $ to get
\begin{equation}
\begin{split}
\iprod{J\,\statevec U^s_t,\statevec U^s}_N = &+ \frac{1}{2} \iprod{\mmatrix S^{-1}\bigmatrix{M}^T\,\mmatrix S\,\spacevec {\mmatrix A}^s\statevec U^s,\spacevec{\nabla}_\xi\,\statevec U^s}_N  - \frac{1}{2} \iprod{\mmatrix S^{-1}\bigmatrix{M}\,\mathrm{S}\,\spacevec{\nabla}_\xi\statevec U^s,(\spacevec {\mmatrix A}^s)^T \statevec U^s}_N \\ 
&- \int\limits_{\partial E,N} (\statevec U^s)^T\left\{\statevec{F}_n^{s,*} -  \frac{1}{2}\statevec{F}^s_n\right\}\,\dS,
\end{split}
\end{equation}
where we define the symmetrized discrete flux $\statevec{F}^s_n$ that uses the symmetric coefficient matrices ${\spacevec {\mmatrix A}^s = \mmatrix S^{-1} \spacevec {\mmatrix A}\;\mmatrix S}$. Using the fact that the symmetrizer matrix $\mmatrix S$ commutes with the metric block matrix $\bigmatrix{M}$  (see e.g. \cite{Gassner_BR1,winters2020construction}) we see that the volume terms cancel out, leaving only surface terms,
\begin{equation}
\begin{split}
\iprod{J\,\statevec U^s_t,\statevec U^s}_N =  &- \int\limits_{\partial E,N} (\statevec U^s)^T\left\{\statevec{F}_n^{s,*} -  \frac{1}{2}\statevec{F}^s_n\right\}\,\dS.
\end{split}
\end{equation}

When we sum over all elements, inner surface terms appear twice (with different normal vectors), whereas element surfaces that are at the physical domain boundary appear only once and are denoted as physical boundary terms ($\PBT$). The interior element surface contributions split into two parts: Surfaces that fall on the material interface $\Gamma$, and those across which the coefficient matrices the same, which we call smooth interface boundary terms, $\operatorname{SIBT}$. The sum over all elements can then be written as
\begin{equation}
\begin{split}
\frac{1}{2}\frac{d}{dt}\sum\limits_{e^k}\inorm{\statevec U^s}^{2}_{J,N}=  &\int\limits_{\Gamma,N} \left\{\jump{ (\statevec U^s)^T}\,\statevec{F}_n^{s,*} - \oneHalf\jump{ (\statevec U^s)^T\statevec{F}^s_n}\right\} \,\dS + \PBT + \operatorname{SIBT},
\end{split}
\end{equation}
written in terms of the jump operator, $\jump{\statevec U} = \statevec U_{R} - \statevec U_{L}$.
Assuming that the discrete physical boundary terms are dissipative, the discrete $L_{2}$ norm satisfies
\begin{equation}
\begin{split}
\frac{1}{2}\frac{d}{dt}\sum\limits_{e^k}\inorm{\statevec U^s}^{2}_{J,N}\leq &\int\limits_{\Gamma,N} \left\{\jump{ (\statevec U^s)^T}\,\statevec{F}_n^{s,*} - \oneHalf\jump{ (\statevec U^s)^T\statevec{F}^s_n}\right\} \,\dS +  \operatorname{SIBT}, 
\end{split}
\end{equation}
which mimics the continuous stability \eqref{eq:GeneralWell-posednessCondition} if $\operatorname{SIBT}\le 0$. 

We thus need a proper numerical flux function $\statevec{F}_n^{s,*}$ to control discrete stability, i.e. to guarantee that the integrand satisfies
\begin{equation}
Q_{N} \equiv\jump{ (\statevec U^s)^T}\,\statevec{F}_n^{s,*} - \oneHalf\jump{ (\statevec U^s)^T\statevec{F}^s_n} \leq Q(\statevec W_{L},\statevec Q_{R})
\label{eq:IntegrandBoundGeneral}
\end{equation}
 pointwise at each node on element faces along the discretization of $\Gamma$, and $\operatorname{SIBT}\le 0$. 

The dissipativity of the $\operatorname{SIBT}$ for the upwind numerical flux has been shown elsewhere, e.g. \cite{Gassner:2013ij},\cite{winters2020construction}. Therefore, in the following we will assume $\operatorname{SIBT}\le 0$ and concern ourselves only with the discontinuous interface terms.

\subsection{Stability for the Scalar Problem}
In the DG approximation, the Rankine-Hugoniot condition and the inflow boundary condition are enforced weakly with the upwind numerical flux
\begin{equation}
F^{*}(U_{L},U_{R}; a_{L},a_{R}) = a_{L} U_{L},
\label{eq:UpwindScalarNumericalFlux}
\end{equation}
 If summation by parts is applied again to the second term in \eqref{eq:DGSEM}, one gets the strong form of the approximation,
in which the integrand of the boundary term is \cite{winters2020construction}
\begin{equation}
F^{*}-F = a_{L}U_{L}- a_{R}U_{R}.
\end{equation}
As the solution converges, this difference goes to zero, and the Rankine-Hugoniot condition is satisfied. Furthermore, 
\begin{equation}
F^{*}-F = a_{L}U_{L}- a_{R}U_{R} = a_{R}\left( \frac{a_{L}}{a_{R}}U_{L} -U_{R} \right),
\end{equation}
so that when the approximation converges, the analytical inflow boundary condition, $U_{R} = \frac{a_{L}}{a_{R}}U_{L}$ is approached, as required, c.f. \eqref{eq:RHScalar}.  

With \eqref{eq:UpwindScalarNumericalFlux}, the interface contribution for the scalar problem is
\begin{equation}
\begin{split}
Q_{N}
& = (U_{R}-U_{L})a_{L} U_{L} - \frac{1}{2}(a_{R}U^{2}_{R}- a_{L}U_{L}^{2})\\&= U_{R}a_{L}U_{L} - a_{L}U_{L}^{2}-\frac{1}{2}a_{R}U^{2}_{R} + \frac{1}{2}a_{L}U_{L}^{2}\\&
= -\frac{1}{2}\left( a_{L}U_{L}^{2} - 2 U_{R}a_{L}U_{L} + a_{R}U^{2}_{R} \right).
\end{split}
\end{equation}
Factoring the quadratic,
\begin{equation}
\begin{split}
Q_{N}&= -\oneHalf\left( a_{L}U_{L}^{2} - 2 U_{R}a_{L}U_{L} + a_{R}U^{2}_{R} \right) 
\\&= 
-\oneHalf a_{L}U_{L}^{2}\left( 1 - 2\frac{U_{R}}{U_{L}} + \frac{a_{R}}{a_{L}}\left(\frac{U_{R}}{U_{L}}\right)^{2}\right)
\\&
= -\oneHalf a_{L}U_{L}^{2}\tilde Q\left(\frac{U_{R}}{U_{L}};a_{L},a_{R} \right).
\end{split}
\end{equation}

The quadratic $\tilde Q(\eta;a_{L},a_{R})$ is concave up and has a minimum when $\eta^{*} = a_{L}/a_{R}$, since
\begin{equation}
\tilde Q'=-2 + 2\frac{a_{R}}{a_{L}} \eta,\quad \tilde Q'' = 2\frac{a_{R}}{a_{L}} >0.
\end{equation}
When $\eta^{*} = a_{L}/a_{R}$, the Rankine-Hugoniot condition is satisfied by the states on either side. The value of that minimum is $Q(\eta^{*}) = 1 - \frac{a_{L}}{a_{R}} $.

It then follows that the contribution to the energy in the numerical approximation matches that of the PDE, \eqref{eq:ScalarEnergyDynamics}, plus a dissipation term dependent on how much the Rankine-Hugoniot condition is not satisfied by the approximate solution. If we define $\beta = a_{L}/a_{R}$, and note that the minimum value of $\tilde Q$ is $1-\beta$, we can separate out that term giving
\begin{equation}
\begin{split}
\tilde Q(\eta;\beta) = 1-2\eta+\frac{1}{\beta}\eta^{2} &= (1-\beta) + (1-2\eta+\frac{1}{\beta}\eta^{2}) - (1-\beta)
\\&
= (1-\beta) + \frac{1}{\beta}(\eta-\beta)^{2}.
\end{split}
\label{eq:QSplittingPDETerm}
\end{equation}
Re-writing the interface contribution in the final form of \eqref{eq:QSplittingPDETerm} will be a key step in showing inferred stability of the approximation for the more complex case of a system of equations.

Then when we substitute for $\eta$ and $\beta$,
\begin{equation}
\begin{split}
Q_{N} &= -\oneHalf a_{L}U_{L}^{2}\left( 1 - \frac{a_{L}}{a_{R}}\right) - \frac{a_{R}U_{L}^{2}}{2}\left( \frac{U_{R}}{U_{L}} -  \frac{a_{L}}{a_{R}}\right)^{2}
\\&
=  -\oneHalf a_{L}\left( 1 - \frac{a_{L}}{a_{R}}\right)U_{L}^{2} -\frac{1}{2a_{R}}\left(a_{R}U_{R}- a_{L}U_{L}\right)^{2}.
\end{split}
\end{equation}

Let us compare:
In the continuous case, we have \eqref{eq:ScalarEnergyDynamics}, with
\begin{equation}
Q\left(u(0^{-}),u(0^{+})\right) = -\frac{a_{L}}{2}\left\{1 -  \frac{a_{L}}{a_{R}}\right\} u^{2}(0^{-}),
\label{eq:QScalarAnalyticRedux}
\end{equation}
whereas discretely,
\begin{equation}Q_{N}=Q\left(U_{L},U_{R}\right) -\frac{1}{2a_{R}}\left(a_{R}U_{R}- a_{L}U_{L}\right)^{2}\le Q\left(U_{L},U_{R}\right).
\label{eq:QScalarDiscreteRedux}
\end{equation}
Thus, according to the definition of stability, Definition \ref{def:StabilityDefinition}, the DGSEM approximation of the scalar problem with the upwind numerical flux has inferred stability.
\begin{remark}
The comparison between \eqref{eq:QScalarAnalyticRedux} and \eqref{eq:QScalarDiscreteRedux} shows explicitly what is interpreted as stability. The first term in \eqref{eq:QScalarDiscreteRedux} can be positive or negative  depending on $a_{L}/a_{R}$, but matches that of the PDE, \eqref{eq:QScalarAnalyticRedux}. The approximation is therefore directly stable if $a_{L}/a_{R}\le 1$, just like the PDE. The second term is always non-positive and represents dissipation of the energy by the approximation.
\qed\end{remark}
\begin{remark}
For the scalar problem it is straightforward to show energy boundedness in a discounted norm by scaling the downwind domain contributions before summing over the elements. When the global sum (in this case, over two elements) is formed,
\begin{equation}
\tilde Q_{\alpha_{c}}=\left( 1 - 2\alpha_{c}\frac{U_{R}}{U_{L}} + \alpha_{c}\frac{a_{R}}{a_{L}}\left(\frac{U_{R}}{U_{L}}\right)^{2}\right).
\end{equation}
As before, $\tilde Q_{N,\alpha_{c}}$ is concave up, with minimum at the same point, $x^{*}$, with minimum value
\begin{equation}
\tilde Q_{min}(a_{L}/a_{R}) = 1 - \alpha_{c}\frac{a_{L}}{a_{R}},
\end{equation}
so 
\begin{equation}
\tilde Q_{\alpha_{c}}(a_{L}/a_{R}) \ge 1 - \alpha_{c}\frac{a_{L}}{a_{R}}.
\label{eq:QTildeLowerBound}
\end{equation}
Since one can always show bounded energy in the new discounted norm by choosing $\alpha_{c}$ to match the analytical value for any (positive) wavespeeds, the condition \eqref{eq:QScalarDiscreteRedux} infers stability. The amount of numerical dissipation in that norm depends on the particular choice of $\alpha_{c}$, however.
\qed\end{remark}
\subsection{Stability for the One-Dimensional Symmetric System}
We now  parallel Sec. \ref{Sec:Symmetric12DsystemPDE} and extend the analysis to a symmetric PDE system in one space dimension. For the system, the DG approximation has the interface contribution
\begin{equation}
Q_{N} = \jump{\statevec U^{T}}\statevec F^{*} - \oneHalf\jump{\statevec U^{T}\mmatrix A\statevec U}.
\label{eq:QNSymmetric}
\end{equation}

The upwind numerical flux is now
\begin{equation}
\begin{split}
\statevec F^{*} &= \mmatrix A_{L}\mmatrix P_{L}\left[\begin{array}{c}\statevec W_{L}^{+}\\\statevec W_{*}^{-}\end{array}\right] = \mmatrix A_{R}\mmatrix P_{R}\left[\begin{array}{c}\statevec W_{*}^{+}\\\statevec W_{R}^{-}\end{array}\right] 
\\&
= \mmatrix P_{L}\Lambda_{L}\left[\begin{array}{c}\statevec W_{L}^{+}\\\statevec W_{*}^{-}\end{array}\right] = 
\mmatrix P_{R}\Lambda_{R}\left[\begin{array}{c}\statevec W_{*}^{+}\\\statevec W_{R}^{-}\end{array}\right],
\end{split}
\end{equation}
with the equalities between the left and right representations arising by virtue of the Rankine-Hugoniot condition. 
Then the key observation is that
\begin{equation}
 \jump{\statevec U^{T}}\statevec F^{*} =  
\statevec U^{T}_{R}\mmatrix P_{R}\Lambda_{R}\left[\begin{array}{c}\statevec W_{*}^{+}\\\statevec W_{R}^{-}\end{array}\right]
- \statevec U^{T}_{L}\mmatrix P_{L}\Lambda_{L}\left[\begin{array}{c}\statevec W_{L}^{+}\\\statevec W_{*}^{-}\end{array}\right].
\end{equation}
But $\statevec U^{T} = \left(\mmatrix P \statevec W\right)^{T} = \statevec W^{T}\mmatrix P ^{T}$ and for the symmetric system $\mmatrix P ^{T}\mmatrix P = \mmatrix I$, so
\begin{equation}
 \jump{\statevec U^{T}}\statevec F^{*} =  
\statevec W^{T}_{R}\Lambda_{R}\left[\begin{array}{c}\statevec W_{*}^{+}\\\statevec W_{R}^{-}\end{array}\right]
- \statevec W^{T}_{L}\Lambda_{L}\left[\begin{array}{c}\statevec W_{L}^{+}\\\statevec W_{*}^{-}\end{array}\right].
\label{eq:jumpUFStar}
\end{equation}
Now,
\begin{equation}
\statevec W^{T}_{L}\Lambda_{L}\left[\begin{array}{c}\statevec W_{L}^{+}\\\statevec W_{*}^{-}\end{array}\right] = 
\statevec W^{+,T}_{L} \bar\Lambda^{+}_{L}\statevec W_{L}^{+} - \statevec W^{-,T}_{L} \left|\bar\Lambda^{-}_{L}\right|\statevec W^{-}_{*}
\end{equation}
and
\begin{equation}
\statevec W^{T}_{R}\Lambda_{R}\left[\begin{array}{c}\statevec W_{*}^{+}\\\statevec W_{R}^{-}\end{array}\right]=
\statevec W^{+,T}_{R} \bar\Lambda^{+}_{R}\statevec W_{*}^{+} - \statevec W^{-,T}_{R} \left|\bar\Lambda^{-}_{R}\right|\statevec W^{-}_{R}.
\end{equation}
Therefore,
\begin{equation}
\jump{\statevec U^{T}}\statevec F^{*} =  
\statevec W^{+,T}_{R} \bar\Lambda^{+}_{R}\statevec W_{*}^{+} 
- \statevec W^{-,T}_{R} \left|\bar\Lambda^{-}_{R}\right|\statevec W^{-}_{R} 
- \statevec W^{+,T}_{L} \bar\Lambda^{+}_{L}\statevec W_{L}^{+} 
+ \statevec W^{-,T}_{L} \left|\bar\Lambda^{-}_{L}\right|\statevec W^{-}_{*}.
\label{eq:jumpUFStarCharacteristic}
\end{equation}

Looking at the second jump term in \eqref{eq:QNSymmetric},
\begin{equation}
\statevec U^{T}\mmatrix A\statevec U = (P\statevec W)^{T}P\Lambda \statevec W = \statevec W^{T} \Lambda \statevec W = \statevec W^{+,T}\bar\Lambda^{+} \statevec W^{+} + \statevec W^{-,T}\bar\Lambda^{-} \statevec W^{-},
\end{equation}
so
\begin{equation}
\begin{split}
\jump{\statevec U^{T}\mmatrix A\statevec U} &= 
\statevec W_{R}^{+,T}\bar\Lambda_{R}^{+} \statevec W_{R}^{+} + \statevec W_{R}^{-,T}\bar\Lambda_{R}^{-} \statevec W_{R}^{-} - 
\statevec W_{L}^{+,T}\bar\Lambda_{L}^{+} \statevec W_{L}^{+} - \statevec W_{L}^{-,T}\bar\Lambda_{L}^{-} \statevec W_{L}^{-}
\\&
= 
\left\{\statevec W_{R}^{+,T}\bar\Lambda_{R}^{+} \statevec W_{R}^{+} 
- \statevec W_{L}^{+,T}\bar\Lambda_{L}^{+} \statevec W_{L}^{+}\right\}
-\left\{ \statevec W_{R}^{-,T}\left|\bar\Lambda_{R}^{-}\right| \statevec W_{R}^{-}
 - \statevec W_{L}^{-,T}\left|\bar\Lambda_{L}^{-}\right| \statevec W_{L}^{-} \right\}.
\end{split}
\end{equation}

Therefore, forming $Q_{N}$ and gathering right and left going wave components,
\begin{equation}\begin{split}
Q_{N}= 
&\left\{
 \statevec W^{+,T}_{R} \bar\Lambda^{+}_{R}\statevec W_{*}^{+} 
- \statevec W^{+,T}_{L} \bar\Lambda^{+}_{L}\statevec W_{L}^{+} 
- \oneHalf \statevec W_{R}^{+,T}\bar\Lambda_{R}^{+} \statevec W_{R}^{+}   
+\oneHalf \statevec W_{L}^{+,T}\bar\Lambda_{L}^{+} \statevec W_{L}^{+}
\right\} 
+
\\&\left\{
\statevec W^{-,T}_{L} \left|\bar\Lambda^{-}_{L}\right|\statevec W^{-}_{*} - \statevec W^{-,T}_{R} \left|\bar\Lambda^{-}_{R}\right|\statevec W^{-}_{R}+ \oneHalf \statevec W_{R}^{-,T}\left|\bar\Lambda_{R}^{-}\right| \statevec W_{R}^{-}
- \oneHalf  \statevec W_{L}^{-,T}\left|\bar\Lambda_{L}^{-}\right| \statevec W_{L}^{-}
\right\}.
\end{split}
\end{equation}
Terms cancel, leaving
\begin{equation}\begin{split}
Q_{N}= 
&-\oneHalf\left\{
 \statevec W^{+,T}_{L} \bar\Lambda^{+}_{L}\statevec W_{L}^{+} 
-2 \statevec W^{+,T}_{R} \bar\Lambda^{+}_{R}\statevec W_{*}^{+} 
+ \oneHalf \statevec W_{R}^{+,T}\bar\Lambda_{R}^{+} \statevec W_{R}^{+}   
\right\} 
\\&
-\oneHalf\left\{
\statevec W_{R}^{-,T}\left|\bar\Lambda_{R}^{-}\right| \statevec W_{R}^{-}
-2\statevec W^{-,T}_{L} \left|\bar\Lambda^{-}_{L}\right|\statevec W^{-}_{*} 
+\oneHalf  \statevec W_{L}^{-,T}\left|\bar\Lambda_{L}^{-}\right| \statevec W_{L}^{-}
\right\}.
\end{split}
\end{equation}
Following \eqref{eq:QSplittingPDETerm}, we now add and subtract terms to match the PDE form, which is
\begin{equation}
\begin{split}
Q 
 = -\oneHalf\left\{\statevec w^{+,T}_{L} \bar\Lambda^{+}_{L}\statevec w^{+}_{L} - \statevec w^{+,T}_{*} \bar\Lambda^{+}_{R}\statevec w^{+}_{*}\right\}
 - \oneHalf\left\{  \statevec w^{-,T}_{R}\left|\bar\Lambda^{-}_{R}\right|\statevec w^{-}_{R} - \statevec w^{-,T}_{*} \left|\bar\Lambda^{-}_{L}\right|\statevec w^{-}_{*}\right\},
\end{split}
\label{eq:QForSystemRepeat}
\end{equation}
to write
\begin{equation}
\begin{split}
Q_{N} = &-\oneHalf
\left\{
\statevec W^{+,T}_{L} \bar\Lambda^{+}_{L}\statevec W^{+}_{L} 
- \statevec W^{+,T}_{*} \bar\Lambda^{+}_{R}\statevec W^{+}_{*}\right\} 
-\oneHalf\left\{
 \statevec W^{+,T}_{*} \bar\Lambda^{+}_{R}\statevec W^{+}_{*} 
-2 \statevec W^{+,T}_{R} \bar\Lambda^{+}_{R}\statevec W_{*}^{+} 
+ \statevec W_{R}^{+,T}\bar\Lambda_{R}^{+} \statevec W_{R}^{+}   
\right\} 
\\&
 - \oneHalf\left\{  \statevec W^{-,T}_{R}\left|\bar\Lambda^{-}_{R}\right|\statevec W^{-}_{R}
  - \statevec W^{-,T}_{*} \left|\bar\Lambda^{-}_{L}\right|\statevec W^{-}_{*}\right\}
\\& -\oneHalf\left\{
\statevec W^{-,T}_{*} \left|\bar\Lambda^{-}_{L}\right|\statevec W^{-}_{*}
-2\statevec W^{-,T}_{L} \left|\bar\Lambda^{-}_{L}\right|\statevec W^{-}_{*} 
+  \statevec W_{L}^{-,T}\left|\bar\Lambda_{L}^{-}\right| \statevec W_{L}^{-}
\right\}.
\end{split}
\end{equation}

Now, let
\begin{equation}
\begin{split}
R^{+} &= \left\{
 \statevec W^{+,T}_{*} \bar\Lambda^{+}_{R}\statevec W^{+}_{*} 
-2 \statevec W^{+,T}_{R} \bar\Lambda^{+}_{R}\statevec W_{*}^{+} 
+ \statevec W_{R}^{+,T}\bar\Lambda_{R}^{+} \statevec W_{R}^{+}   
\right\} ,
\\
R^{-} &= \left\{
\statevec W^{-,T}_{*} \left|\bar\Lambda^{-}_{L}\right|\statevec W^{-}_{*}
-2\statevec W^{-,T}_{L} \left|\bar\Lambda^{-}_{L}\right|\statevec W^{-}_{*} 
+  \statevec W_{L}^{-,T}\left|\bar\Lambda_{L}^{-}\right| \statevec W_{L}^{-}
\right\}.
\end{split}
\end{equation}
Then
\begin{equation}
\begin{split}
Q_{N} = &-\oneHalf
\left\{
\statevec W^{+,T}_{L} \bar\Lambda^{+}_{L}\statevec W^{+}_{L} 
- \statevec W^{+,T}_{*} \bar\Lambda^{+}_{R}\statevec W^{+}_{*}\right\} - \oneHalf R^{+}
\\&
- \oneHalf\left\{  \statevec W^{-,T}_{R}\left|\bar\Lambda^{-}_{R}\right|\statevec W^{-}_{R}
  - \statevec W^{-,T}_{*} \left|\bar\Lambda^{-}_{L}\right|\statevec W^{-}_{*}\right\} -\oneHalf R^{-}.
\end{split}
\end{equation}

To show that the approximation is stable according to Definition \ref{def:StabilityDefinition}, then, we just need to show that $R^{\pm}\ge 0$, since the other terms match those of the PDE.
To do so, let $\bar {\statevec W}^{\pm} = \sqrt{\left| \bar \Lambda^{\pm}\right|}\statevec W^{\pm}$.
Then
\begin{equation}
R^{+}= \bar{\statevec W}^{+,T}_{*} \bar{\statevec W}^{+}_{*} 
-2 \bar{\statevec W}^{+,T}_{R} \bar{\statevec W}_{*}^{+} 
+ \bar{\statevec W}_{R}^{+,T}\bar{ \statevec W}_{R}^{+} = \left(  \bar{\statevec W}^{+}_{*} - \bar{ \statevec W}_{R}^{+} \right)^{2}\ge 0.
\end{equation}
Similarly,
\begin{equation}
R^{-}= \left(  \bar{\statevec W}^{-}_{*} - \bar{ \statevec W}_{L}^{-} \right)^{2}\ge 0.
\end{equation}
Thus, the interface contribution matches that of the PDE plus an additional dissipation and has inferred stability, satisfying Definition \ref{def:StabilityDefinition}
with
\begin{equation}
\resizebox{.9\hsize}{!}{$
Q_{N} \le -\oneHalf
\left\{
\statevec W^{+,T}_{L} \bar\Lambda^{+}_{L}\statevec W^{+}_{L} 
- \statevec W^{+,T}_{*} \bar\Lambda^{+}_{R}\statevec W^{+}_{*}\right\}
- \oneHalf\left\{  \statevec W^{-,T}_{R}\left|\bar\Lambda^{-}_{R}\right|\statevec W^{-}_{R}
  - \statevec W^{-,T}_{*} \left|\bar\Lambda^{-}_{L}\right|\statevec W^{-}_{*}\right\} = Q(\statevec W_{L}, \statevec W_{R}).$}
\label{eq:QNFor1DSystemBound}
\end{equation}

\subsection{Stability of the General Problem}
As in the continuous problem, we use the analysis of the one-dimensional problem to imply stability of the multidimensional one. As before, replace $\statevec U \leftarrow \statevec U^{s}$ and $\mmatrix A\leftarrow \spacevec {\mmatrix A}^{s}\cdot\spacevec n$. Then $Q_{N}$ is given by \eqref{eq:QNFor1DSystemBound}, with the eigenvalues (and eigenvectors to construct the characteristic variables) coming from $\spacevec {\mmatrix A}^{s}\cdot\spacevec n$. Therefore the approximation to the general multidimensional problem is stable according to Definition \ref{def:StabilityDefinition}.

\section{Example}

As an example, we consider the scattering of a plane wave off a plane material interface, approximating the system of equations \eqref{eq:DivergenceEquation} with the state vector and coefficient matrices \eqref{eq:WaveEqnMatrices} reduced to two space dimensions.
The problem has exact incident, transmitted and reflected plane wave solutions of the form
\begin{equation}
\statevec u = a\psi\left(\spacevec k\cdot\spacevec x - \omega(t - t_{0})\right)\left[\begin{array}{c}1 \\\frac{k_x}{\rho c} \\\frac{k_y}{\rho c}\end{array}\right],
\label{eq:WavePacketSolution}
\end{equation}
where $\psi$ is a given wavefunction, $a$ is the amplitude, $\spacevec k$ is the wavevector, $\omega$ is the frequency. For the incident wavevector 
\begin{equation}
\spacevec k^i  = \frac{\omega }{{c_L }}\left( {k_x^i \hat x + k_y^i \hat y} \right),
\end{equation}
the reflected and transmitted wavevectors are
\begin{equation}
\begin{array}{l}
 \spacevec k^r  = \dfrac{\omega }{{c_L }}\left( { - k_x^i \hat x + k_y^i \hat y} \right) \\ \\
 \spacevec k^T  = \dfrac{\omega }{{c_R }}\left[ {\sqrt {1 - \left( {\dfrac{{c_R }}{{c_L }}} \right)^2 \left( {k_y^i } \right)^2 } \hat x + \dfrac{{c_R }}{{c_L }}k_y^i \hat y} \right], \\ 
 \end{array}
\end{equation}
with amplitudes
\begin{equation}
\begin{array}{l}
 \dfrac{{a^r }}{{a^i }} = \dfrac{1}{d}\left( {\rho _R c_R k_x^T /|\spacevec k^T|  - \rho _L c_L k_x^i /|\spacevec k^i| } \right), \\ \\
 \dfrac{{a^T }}{{a^i }} = \dfrac{1}{d}\left( {\rho _L c_L k_x^r /|\spacevec k^r|  - \rho _R c_L k_x^i /|\spacevec k^i |} \right), \\ 
 \end{array}
\end{equation}
where
\begin{equation}
d =  - \rho _R c_R k_x^T /|\spacevec k^T|  + \rho _L c_L k_x^r /|\spacevec k^r| .
\end{equation}
For the wavefunction, we choose the Gaussian
 \begin{equation}
 \psi(s) = e^{-s^{2}/(\omega\sigma)^{2}},
 \end{equation}
 with $\sigma^{2} = -(MT)^{2}/(4\ln(10^{-4}))$, $M = 4$ and period $T = 2\pi/\omega$.

We compute the problem on the square domain $[-5,5]^{2}$ with 400 square elements and the material interface at $x=0$. The solution parameters are provided in Table \ref{tab:ReflTransTable}.
	 \begin{table}[tbp]
	\begin{center}
	\caption {Parameters for Plane Wave Reflection Problem}
	\label{tab:ReflTransTable}
	\begin{tabular}{c|ccccccccc}
	Parameter & $M$ & $\omega$ & $k^i_x$ & $k^i_y$ & $\rho_L$ & $\rho_R$ & $c_L$ & $c_R$ & $t_{0}$ \\
	\hline
	Value & 4 & $4\pi$ & 0.5 & $\sqrt{3/2}$ & 1 & 0.4 & 1 & 0.7 & 3
	\end{tabular}
	\end{center}
	\end{table}%

\begin{figure}[htbp] 
   \centering
   \includegraphics[width=3in]{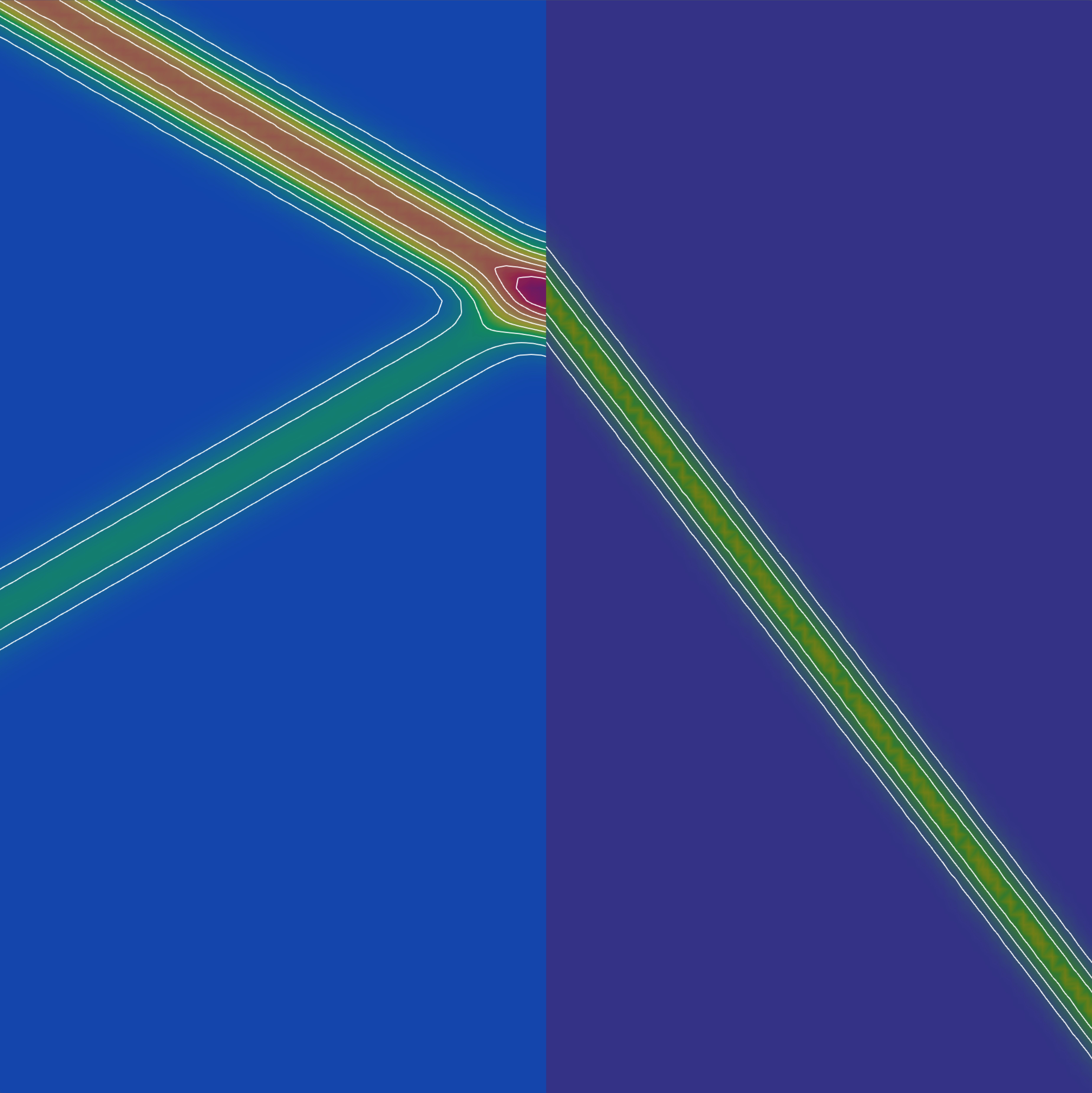} 
   \caption{Computed $p$ contours at time $t=5$ for plane wave scattering from a material interface along the vertical center of the domain}
   \label{fig:pContours}
\end{figure}

\begin{figure}[htbp] 
   \centering
   \includegraphics[width=3in]{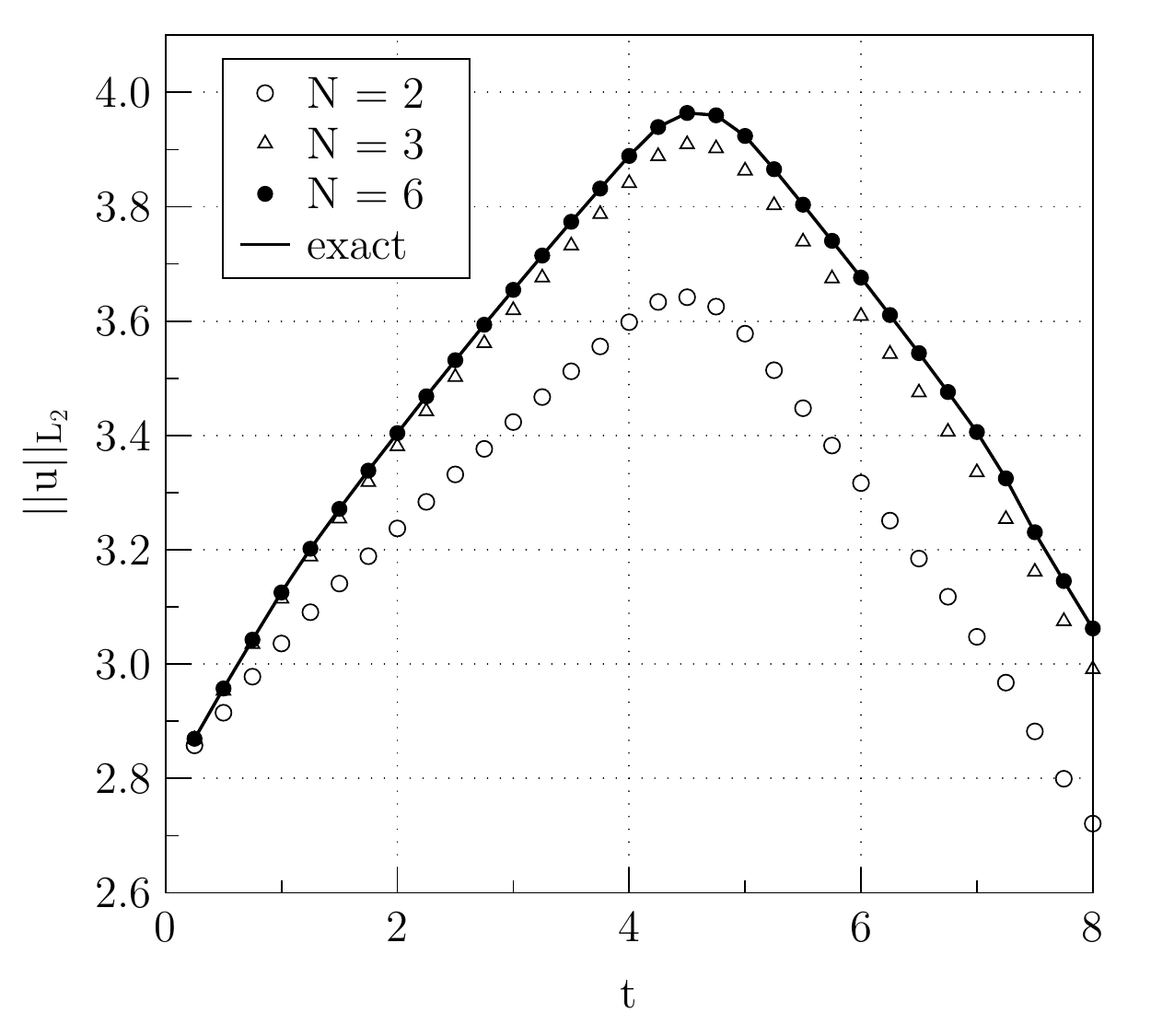} 
   \caption{$L_{2}$ energy as a function of time for scattering at a material interface }
   \label{fig:ComputedAndExactEnergy}
\end{figure}

The results are shown in Figs. \ref{fig:pContours}  and \ref{fig:ComputedAndExactEnergy}. Fig. \ref{fig:pContours} shows the contours of the $p$ component of the solution at time $t=5.0$, which is near the time of the maximum $L_{2}$ energy,  computed with sixth order polynomials. Clearly seen is the jump discontinuity at the interface. The $L_{2}$ energy is plotted as a function of time in Fig. \ref{fig:ComputedAndExactEnergy}, for polynomial degrees $N= 2, 3$ and 6. Although the $L_{2}$ energy initially grows, it reaches a maximum around time $t=4.5$. Fig. \ref{fig:ComputedAndExactEnergy} shows that the computed energy converges from below to the exact as the polynomial order is increased. In fact, it converges exponentially with polynomial degree, as expected \cite{canuto2006} for a spectral element method. Also, as expected due to the additional dissipation at physical, smooth and discontinuous interfaces, the computed energies fall below the exact curve and are worst for low order approximations.

\section{Conclusions}

We have shown that the interface treatment of the discontinuous Galerkin spectral element method with the upwind numerical flux is stable for hyperbolic systems with discontinuous coefficient matrices when the eigenvectors are preserved across the interface. Examples include systems like Maxwell's equations, or acoustic and elastic wave equations. The new feature of our approach was to show that the discrete $L_{2}$ norm of the approximate solution grows no faster than the same norm of the continuous solution. By matching the $L_{2}$ norm, we avoid having to find the precise conditions for a discounted norm in which the energy is bounded by the initial data (for homogenous and dissipative boundary conditions). The numerical flux only weakly enforces the inflow boundary condition and the Rankine-Hugoniot condition. Viewing stability in terms of the $L_{2}$ norm shows that the dissipation introduced by the upwind numerical flux depends on the amount by which the approximate solution fails to satisfy the Rankine-Hugoniot condition.

\acknowledgement
{
The authors would like to thank Andrew Winters, Lucas Wilcox and Juan Manzanero for helpful advice. 
 This work was supported by a grant from the Simons Foundation (\#426393, David Kopriva). Gregor Gassner thanks the Klaus-Tschira Stiftung and
the European Research Council for funding through the ERC Starting Grant “An
Exascale aware and Un-crashable Space-Time-Adaptive Discontinuous Spectral
Element Solver for Non-Linear Conservation Laws” (EXTREME, project no. 71448). Jan Nordstr\"om was supported by Vetenskapsrådet, Sweden grant nr: 2018-05084 VR and the Swedish e-Science Research Center (SeRC).
}
\bibliographystyle{plain}
\bibliography{jumpPaper.bib}

\begin{thebibliography}{10}

\bibitem{Isi:A1981Lw20700001}
S~Abarbanel and D~Gottlieb.
\newblock {Optimal Time Splitting For Two-Dimensional And 3-Dimensional
  {N}avier-{S}tokes Equations With Mixed Derivatives}.
\newblock {\em {Journal Of Computational Physics}}, {41}({1}):{1--33}, {1981}.

\bibitem{canuto2006}
C.~Canuto, M.~Hussaini, A.~Quarteroni, and T.~Zang.
\newblock {\em Spectral Methods: Fundamentals in Single Domains}.
\newblock Springer, Berlin, 2006.

\bibitem{ISI:000226090600009}
S.Z. Deng, W~Cai, and V.N. Astratov.
\newblock {Numerical study of light propagation via whispering gallery modes in
  microcylinder coupled resonator optical waveguides}.
\newblock {\em OPTICS EXPRESS}, {12}({26}):{6468--6480}, {DEC 27} {2004}.

\bibitem{Gassner_BR1}
G.~J. Gassner, A.~R. Winters, F.~J. Hindenlang, and D.~A. Kopriva.
\newblock The {BR1} scheme is stable for the compressible {N}avier-{S}tokes
  equations.
\newblock {\em {Journal of Scientific Computing}}, 77(1):154--200, 2018.

\bibitem{doi:10.1137/16M1087710}
Fatemeh Ghasemi and Jan Nordstr{\"o}m.
\newblock Coupling requirements for multiphysics problems posed on two domains.
\newblock {\em SIAM Journal on Numerical Analysis}, 55(6):2885--2904, 2017.

\bibitem{Hesthaven:2002uq}
J.~S. Hesthaven and T.~Warburton.
\newblock Nodal high-order methods on unstructured grids. {I.} {T}ime-domain
  solution of {M}axwell's equations.
\newblock {\em Journal of Computational Physics}, 181:186--221, 2002.

\bibitem{Kopriva:2006er}
D.~A. Kopriva.
\newblock Metric identities and the discontinuous spectral element method on
  curvilinear meshes.
\newblock {\em The Journal of Scientific Computing}, 26(3):301--327, March
  2006.

\bibitem{Gassner:2013ij}
D.~A. Kopriva and G.~Gassner.
\newblock An energy stable discontinuous {G}alerkin spectral element
  discretization for variable coefficient advection problems.
\newblock {\em {SIAM Journal on Scientific Computing}}, 36(4):A2076--A2099,
  2014.

\bibitem{Koprivaetal1999}
D.~A. Kopriva, S.~L. Woodruff, and M.~Y. Hussaini.
\newblock Discontinuous spectral element approximation of {M}axwell's
  {E}quations.
\newblock In B.~Cockburn, G.~Karniadakis, and C.-W. Shu, editors, {\em
  Proceedings of the {I}nternational {S}ymposium on {D}iscontinuous {G}alerkin
  {M}ethods}, pages 355--361, New York, May 2000. Springer-Verlag.

\bibitem{Kopriva:2009nx}
David~A. Kopriva.
\newblock {\em Implementing Spectral Methods for Partial Differential
  Equations}.
\newblock Scientific Computation. Springer, May 2009.

\bibitem{10.1007/978-3-319-65870-4_2}
David~A. Kopriva.
\newblock A polynomial spectral calculus for analysis of {DG} spectral element
  methods.
\newblock In Marco~L. Bittencourt, Ney~A. Dumont, and Jan~S. Hesthaven,
  editors, {\em Spectral and High Order Methods for Partial Differential
  Equations ICOSAHOM 2016}, pages 21--40, Cham, 2017. Springer International
  Publishing.

\bibitem{Kopriva2016274}
David~A. Kopriva and Gregor~J. Gassner.
\newblock Geometry effects in nodal discontinuous {G}alerkin methods on curved
  elements that are provably stable.
\newblock {\em Applied Mathematics and Computation}, 272, Part 2:274 -- 290,
  2016.

\bibitem{La-Cognata:2016ng}
Cristina La~Cognata and Jan Nordstr{\"o}m.
\newblock Well-posedness, stability and conservation for a discontinuous
  interface problem.
\newblock {\em BIT Numerical Mathematics}, 56(2):681--704, 2016.

\bibitem{Manzanero:2018kk}
Juan Manzanero, Gonzalo Rubio, Esteban Ferrer, Eusebio Valero, and David~A.
  Kopriva.
\newblock Insights on aliasing driven instabilities for advection equations
  with application to {G}auss--{L}obatto discontinuous {G}alerkin methods.
\newblock {\em Journal of Scientific Computing}, 75(3):1262--1281, 2018.

\bibitem{Nordstrom:2016jk}
Jan Nordstr\"om.
\newblock A roadmap to well posed and stable problems in computational physics.
\newblock {\em {Journal Of Scientific Computing}}, DOI
  10.1007/s10915-016-0303-9, 2016.

\bibitem{wilcox2010}
L.~C. Wilcox, G.~Stadler, C.~Burstedde, and O.~Ghattas.
\newblock A high-order discontinuous {G}alerkin method for wave propagation
  through coupled elastic-acoustic media.
\newblock {\em Journal of Computational Physics}, 229(24):9373--9396, 2010.

\bibitem{Winters:2013nx}
Andrew~R. Winters and David~A. Kopriva.
\newblock {ALE-DGSEM} approximation of plane wave reflection and transmission
  from a moving medium.
\newblock {\em {Journal of Computational Physics}}, 263(1):176--202, 2014.

\bibitem{winters2020construction}
Andrew~R. Winters, David~A. Kopriva, Gregor~J. Gassner, and Florian Hindenlang.
\newblock Construction of modern robust nodal discontinuous {G}alerkin spectral
  element methods for the compressible {N}avier-{S}tokes equations, 2020.

\end{thebibliography}

\end{document}